\documentclass[12pt]{article}
\usepackage{amsfonts}
\usepackage{amssymb}
\usepackage{polski}

\input{tcilatex}
\begin{document}

\begin{center}
\textbf{Set Theory and the Analyst}

by

\textbf{N. H. Bingham and A. J. Ostaszewski}

\bigskip

\bigskip

{\small \textit{Then to the rolling heaven itself I cried},}\\[0pt]
{\small \textit{Asking what lamp had destiny to guide} }\\[0pt]
{\small \textit{Her little children stumbling in the dark}.}\\[0pt]
{\small \textit{And `A blind understanding' heaven replied.}\\[0pt]
-- The Rubaiyat of Omar Khayyam }

\bigskip
\end{center}

\textbf{Abstract. }This survey is motivated by specific questions arising in
the similarities and contrasts between (Baire) category and (Lebesgue)
measure -- category-measure duality and non-duality, as it were. The bulk of
the text is devoted to a summary, intended for the working analyst, of the
extensive background in set theory and logic needed to discuss such matters:
to quote from the Preface of Kelley [Kel]: "what every young analyst should
know".

\bigskip

\begin{center}
\textbf{Table of Contents}
\end{center}

\noindent 1. Introduction\newline
\noindent 2. Early history\newline
\noindent 3. G\"{o}del Tarski and their legacy\newline
\noindent 4. Ramsey, Erd\H{o}s and their legacy: infinite combinatorics%
\newline
\qquad \hbox{          }\hbox{          }4a. Ramsey and Erd\H{o}s\newline
\qquad \hbox{          }\hbox{          }4b. Partition calculus and large
cardinals\newline
\qquad \hbox{          }\hbox{          }4c. Partitions from large cardinals%
\newline
\qquad \hbox{          }\hbox{          }4d. Large cardinals continued

\noindent 5. Beyond the constructible hierarchy $L$ -- I\newline
\qquad \hbox{          }\hbox{          }5a. Expansions via ultrapowers%
\newline
\qquad \hbox{          }\hbox{          }5b. Ehrenfeucht-Mostowski models:
expansion via indiscernibles\newline
\noindent 6. Beyond the constructible hierarchy $L\ $-- II\newline
\qquad \hbox{          }\hbox{          }6a. Forcing and generic extensions%
\newline
\qquad \hbox{          }\hbox{          }6b. Forcing Axioms

\noindent 7. Suslin, Luzin, Sierpi\'{n}ski and their legacy: infinite games
and large cardinals\newline
\qquad \hbox{          }\hbox{          }7a. Analytic sets.\newline
\qquad \hbox{          }\hbox{          }7b. Banach-Mazur games and the
Luzin hierarchy\newline
\noindent 8. Shadows\newline
\noindent 9. The syntax of Analysis: Category/measure regularity versus
practicality\newline
\noindent 10. Category-Measure duality\newline
\noindent Coda

\bigskip

\textbf{1. Introduction}

An analyst, as Hardy said, is a mathematician habitually seen in the company
of the real or complex number systems. For simplicity, we restrict ourselves
to the reals here, as the complex case is obtainable from this by a
cartesian product. In mathematics, one's underlying assumption is generally
Zermelo-Fraenkel set theory (ZF), augmented by the Axiom of Choice (AC) when
needed, giving ZFC. One should always be as economical as possible about
one's assumptions; it is often possible to proceed without the full strength
of AC. Our object here is to survey, with the working analyst in mind, a
range of recent work, and of situations in analysis, where one can usefully
assume less than ZFC. This survey is rather in the spirit of that by Wright
[Wri] forty years ago and Mathias's `Surrealist landscape with figures'
survey [Mat]; a great deal has happened since in the area, and we feel that
the time is ripe for a further survey along such lines.

As classical background for what follows, we refer to the book of Oxtoby
[Oxt] on (Lebesgue) measure and (Baire) category. This book explores the
duality between them, focussing on their remarkable similarity; for Oxtoby,
it is the measure case that is primary. Our viewpoint is rather different:
for us, it is the \textit{category} case that is primary, and we focus on
their differences. Our motivation was that a number of results in which
category and measure behave interchangeably \textit{disaggregate} on closer
examination. One obtains results in which what can be said depends
explicitly on what axioms one assumes. One thus needs to adopt a flexible,
or pluralist, approach in order to be able to handle apparently quite
traditional problems within classical analysis.

We come to this survey with the experience of a decade of work on problems
with wide-ranging contexts, from the real line to topological groups, in
which `the category method' has been key. The connection with the Baire
Category Theorem, viewed as an equivalent of the Axiom of Dependent Choices
(DC: for every non-empty set $X$ and $R\subseteq X^{2}$ satisfying $\forall
x\exists y[(x,y)\in R]$ there is $h\in X^{\mathbb{N}}$ with $%
(h(n),h(n+1))\in R$ for all $n\in \mathbb{N})$, has sensitized us to a
reliance on this as a very weak version of the Axiom of Choice, AC, but one
that is often adequate for analysis. This sensitivity has been further
strengthened by settings of general theorems on Banach spaces reducible to
the separable case (e.g. via Blumberg's Dichotomy [BinO6, Th. B]; compare
the separable approximations in group theory [MonZ, Ch. II \S 2.6]) where DC
suffices. On occasion, it has been possible to remove dependence on the
Hahn-Banach theorem, a close relative of AC.

For the relative strengths of the usual Hahn-Banach Theorem HB and the Axiom
of Choice AC, see [Pin1,2]; [PinS] provide a model of set theory in which
the Axiom of Dependent Choices DC holds but HB\ fails. HB\ is derivable from
the Prime Ideal Theorem PI, an axiom weaker than AC : for literature see
again [Pin1,2]; for the relation of the Axiom of Countable Choice, ACC
(below), to DC\ here, see [HowR]. Note that HB for separable normed spaces
is not provable from DC [DodM, Cor. 4], unless the space is complete --- see
[BinO10].

When category methods fail, e.g. on account of `character degradation', as
when the \textit{limsup} operation is applied to well-behaved functions (see 
\S 9), the obstacles may be removed by appeal to supplementary set-theoretic
axioms, so leading either beyond, or sometimes away from, a classical
setting. This calls for analysts to acquire an understanding of their
interplay and their standing in relation to `classical intuition' as
developed through the historical narrative. Our aim here is to describe this
hinterland in a language that analysts may appreciate.

We list some sources that we have found useful, though we have tried to make
the text reasonably self-contained. From logic and foundations of
mathematics, we need AC and its variants, for which we refer to Jech [Jec1].
For set theory, our general needs are served by [Jec2]; see also Ciesielski
[Cie], Shoenfield [Sho], Kunen [Kun3]. For descriptive set theory, see
Kechris [Kec2] or [MarK]; for analytic sets see Rogers et al. [Rog]. For
large cardinals, see e.g. Drake [Dra], Kanamori [Kan], Woodin [Woo1].

The paper is organised as follows. After a review of the early history of
the axiomatic approach to set theory (including also a brief review of some
formalities) we discuss the contributions of G\"{o}del and Tarski and their
legacy, then of Ramsey and Erd\H{o}s and their legacy. We follow this with a
discussion of the role of infinite combinatorics (partition calculus) and of
the `large cardinals'. We then sketch the various `pre-Cohen' expansions of $%
\ $G\"{o}del's universe of constructible sets $L$ (via the ultrapowers of \L %
o\'{s}, or the indiscernibles of Ehrenfeucht-Mostowski models, and the
insights they bring to our understanding of $L$). This is followed by an
introduction to the `forcing method' and the generic extensions which it
enables. We describe classical descriptive set theory: the early
`definability theme' pursued by Suslin, Luzin, Sierpi\'{n}ski and their
legacy, and the completion of their programme more recently through a
recognition of the unifying role of Banach-Mazur games; these require large
cardinals for an analysis of their `consistency strength', and are seen in
some of the most recent literature as casting `shadows' (\S 8) on the real
line.

We conclude in \S 9 with a discussion of the `syntax of analysis' in order
to draw on the `definability theme', and in that light we turn finally in \S %
10 to the additional axioms which permit a satisfying category-measure
duality for the working analyst.

\noindent \textit{The canonical status of the reals}

The English speak of `the elephant in the room', meaning something that
hangs in the air all around us, but is not (or little) spoken of in polite
company. The canonical status of the reals is one such, and the \textit{sets}
of reals available is another such. We speak of `the reals' (or `the real
line'), $\mathbb{R}$ and `the rationals', $\mathbb{Q}$ -- as in the Hardy
saying above. The definite article suggests canonical status, in some sense
-- what sense?

The rationals are indeed canonical. We may think of them as a `tie-rack', on
which irrationals are `hung'. But how, and how many? In \S 6 we will review
the forcing method for selecting from `outside the wardrobe'\ an initial lot
of almost any size that one may wish for, together with their `Skolem hull'
-- further ones required by the operation of the axioms (see the Skolem
functions of \S 2 below).

In brief: the reals are \textit{canonical modulo cardinality, but not
otherwise}. This is no surprise, in view of the Continuum Hypothesis (CH),
which directly addresses the cardinality of the real line/continuum, and
which we know from P. J. Cohen's work of 1963/64 [Coh1,2] will always be
just that, a hypothesis. The canonical status of the reals rests on (at
least) four things:

(i) (geometrical; ancient Greeks): lines in Euclidean geometry: any line can
be made into a cartesian axis;

(ii) (analytical; 1872): Dedekind cuts;

(iii) (analytical; 1872): Cantor, equivalence classes of Cauchy sequences of
rationals (subsequently extended topologically: completion of any metric
space);

(iv) (algebraic; modern): any \textit{complete} archimedean ordered field is
isomorphic to $\mathbb{R}$ (see e.g. [Cohn, \S\ 6.6]: our italics; here too
Cauchy sequences are used to define completeness).

None of these is concerned with cardinality; CH is. On the other hand,
completeness depends on which $\omega $-sequences are available.
Accordingly, the problems that confront the working analyst split, into two
types. Some (usually the `less detailed') do not hinge on cardinality, and
for these the reals retain their traditional canonical status. By contrast,
some do hinge on cardinality; these are the ones that lead the analyst into
set-theoretic underpinnings involving an element of \textit{choice}. Such
choices emphasise the need for a \textit{plural} approach, to axiomatic
assumptions, and hence to the status of the reals. This is inevitable: as
Solovay [Sol1] puts it, `it (the cardinality of the reals) can be anything
it ought to be'.

We turn now to the second of the `elephants' above: which sets of reals are
available. The spectrum of axiom possibilities which we review in \S 10
extends from the `prodigal' (below -- see \S 2) AC at one end (which yields
for example non-measurable Vitali sets) to the restrictive DC\ with
additional components of LM\ (`all sets of reals are measurable') and/or PB
(`all sets of reals have the Baire property') at the other, and include
intermediate positions for the additional component such as PD\ (`all
projective sets of reals are determined'), where the sets of reals with
these so-called `regularity properties' are qualified (see \S \S 7 and 9).

Underlying an analysis of these axioms is repeated appeal to simplification
of contexts - a mathematical \textit{ex oriente lux} -- typified by passage
to a `large' homogeneous/monochromatic subset, as in Ramsey's Theorem on $%
\mathbb{N}$ (\S 4a). This has generalizations to large cardinals $\kappa $,
in particular ones that support a $\{0,1\}$-valued measure (equivalently, a
`suitably complete' ultrafilter -- see below). On the one hand, the latter
permits an extension of Suslin's classical tree-like representation of an
analytic set (\S 7) to sets of far greater logical complexity by witnessing
membership of a set by means of infinite branches in a corresponding tree
that pass through `large' sets of nodes at each height/level (see \S 7). On
the other hand, in the context of the `line' of ordinals, one meets other
forms of isomorphic behaviour on `large' sets: on closed unbounded subsets
of ordinals and on the related stationary sets (\S 5b, 6b).

\textit{Notes.} 1. This survey arose out of our decade-long probing of
questions in regular variation [BinGT]. In [BinO3] we needed to disaggregate
a classical theorem of Delange (see [BinGT, Th. 2.01]); the category and
measure aspects need different set-theoretic assumptions. We regard the
category case as primary, as one can obtain the measure case from it by
working bitopologically (passing from the Euclidean to the \textit{density
topology}; see [BinO2,7,8]); also, measure theory needs stronger
set-theoretic assumptions than category theory (\S 10.2 and \S 10.3 below).
If one replaces the limits in regular variation by \textit{limsups}, the
Baire property or measurability may be lost; the resulting character
degradation is studied in detail in [BinO3 \S 3, 5\S 11].

2. We close by a brief mention of `yet another elephant in the room'. One
can never prove consistency (of sets of rich enough axioms), merely relative
consistency. This is related to G\"{o}del's incompleteness theorems (\S 3).
Thus we do not know that ZF or ZFC itself is consistent; this is something
we have to live with; it is no reason to despair, or give up mathematics;
quite the contrary, if anything. In what follows, `consistency' means
`consistency relative to ZF'.

\bigskip

\textbf{2. Early history}

A little historical background may not come amiss here. The essence of
analysis -- and the reason behind the Hardy quotation that we began with --
is its concern with infinite or limiting processes -- most notably, as in
calculus, our most powerful single technique in mathematics (and indeed, in
science generally). Life being only finitely long, the infinite -- actual or
potential -- takes us beyond direct human experience, even in principle.
This underlies the unease the ancient Greeks had with the irrationals (or
reals), and why they missed calculus (at least in its differential form,
despite their success with areas and volumes under the heading of the
`method of exhaustion'). One can see, for example in the ordering of the
material in the thirteen books of Euclid's \textit{Elements, }that they were
at ease with rationals, and with geometrical objects such as line-segments
etc., but not with reals. Traces of this unease survive in Newton's handling
of the material in his \textit{Principia, }where he was at pains to use
established geometrical arguments rather than his own `method of fluxions'.
That there was unfinished business here shows, e.g., in the title of a work
of one of the founding fathers of analysis, Bolzano, with his \textit{%
Paradoxien des Unendlichen }(1852, posthumous). The bridge between the real
line and the complex plane (the `Argand diagram' -- Argand, 1806, Wessel,
1799, Gauss, 1831) pre-dated this. The construction of the reals came
independently in two different ways in 1872: Dedekind cuts (or sections),
which still dominate settings where one has an \textit{order, }and Cantor's
construction via (equivalence classes of) Cauchy sequences (of rationals) --
still ubiquitous, as the completion procedure for metric spaces.

\textit{Cantor. }Cantor's work, in the 1870s to 1890s, established set
theory (\textit{Mengenlehre)} as the basis on which to do mathematics, and
analysis in particular. Here we find, for example, the countability of the
rationals, and of the algebraic numbers (Cantor, 1874) and the
uncountability of the reals (Cantor, 1895), established via the familiar
Cantor diagonalisation argument. But note what is implicit here: Cantor
diagonalisation (as used, say, to prove the countability of the rationals)
is an \textit{effective }argument. But to move from this to saying that `the
union of countably many countable sets is countable' (Cantor, 1885) needs
the Axiom of Countable Choice (ACC), below.

\textit{Hilbert. }Moving to the 20th century: Hilbert famously said (in
defence of Cantor against Kronecker): `No one shall expel us from the
paradise that Cantor has created for us'. Hilbert addressed himself to the
programme of re-working the mathematical canon of its time to (then) modern
standards of rigour, witness his books on the foundations of geometry
[Hil1,2,3] (1899) and of mathematics [HilB] (1934, 1939, with Bernays), cf.
the Hilbert problems of 1900. As we shall see, Hilbert was a man of his time
here, and his views on foundational questions were too naive. Meanwhile,
Lebesgue introduced measure theory in 1902, Fr\'{e}chet metric spaces in
1906, and Hausdorff general topology in 1905-1914 (three very different
editions of his classic book \textit{Grundzuge der Mengenlehre }appeared in
1914, 1927 and 1935). Hilbert space emerged c. 1916 (work of Hilbert and
Schmidt; named by F. Riesz in 1926). Banach's book [Ban] appeared in 1932,
effectively launching the field of functional analysis; this magisterial
work is still worth reading. But, Banach was a man of his time; he worked
sequentially, rather than using the language of weak topologies, presumably
because he felt it to be not yet in final form. However, the language and
viewpoint of general topology was already available, and already a
speciality of the new Polish school of mathematics, of which Banach himself
was the supreme ornament. For a scholarly and sympathetic account of these
matters, see Rudin [Rud, Appendix B].

The need for care in set theory had been dramatically shown by the Russell
Paradox of 1902, and its role in showing the limitations of Frege's
programme in logic and foundations, especially his \textit{Grundgesetze der
Arithmetik} (vol. 2 of 1903). The Paradox, far from being a programme
wrecker, was pregnant with consequences [GabW], just as with G\"{o}del's
work later (below), and that too was ultimately based on a Paradox (the
`Liar paradox'). See [Hall2] for a discussion. Foundational questions had
been addressed in 1889 by Peano. Zermelo began his axiomatisation, and gave
the Axiom of Choice (AC) in 1902. Fraenkel, Skolem and others continued and
revised this work; what is known nowadays as Zermelo-Fraenkel set theory
(ZF), together with ZF+AC, or ZFC, emerged by 1930 or so. AC is most often
used in the (equivalent) form of Zorn's Lemma of 1935 (a misnomer, as the
result is due to Kuratowski in 1922, but the usage is now established). It
will be helpful for later passages to note that the axioms include the
operations of comprehension (the forming of a subset determined by a
property), union and power set (denoted here by $\wp $), as well as
foundation/regularity, asserting the well-foundedness of the relation of
membership $\in $ (no descending $\in $-chains). In this context AC\ is a
generator of sets par excellence, with effects of both positive and negative
aspects: allowing the construction both to `satisfy intuition' (as in the
construction of `invariant means') and to astound it (as in the
Banach-Tarski paradox): see the comments in [TomW, Ch. 15]. The tension
between `too many' sets or `too few' pervades the history of set theory
through the lens of logic, all the way back to Cantor: see [Hall1]. For a
discussion of approaches to axiomatization see [Sco2].

\textit{Brouwer. }The interplay between analysis (specifically, topology)
and foundations in this era is well exemplified by the work of Brouwer.
Brouwer is best remembered for two contributions: his fixed-point theorem
(of 1911, [Bro1]), and Intuitionism (1920, cf. [Bro2]). The first is beloved
of economists, as it provides existence proofs of economic equilibria -- the
`invisible hand' of Adam Smith, and his later `disciples'. But, his proof of
the fixed-point theorem was a non-constructive existence proof, and Brouwer
lost faith in these for foundational reasons. He reacted by seeking to
re-formulate mathematics `intuitively', on new foundations -- differing from
those in use then and now by, for instance, outlawing proof by
contradiction. This led to serious conflict, for instance the \textit{%
Annalenstreit }(Annals struggle) of 1928, where Hilbert, as Editor-in-Chief
of the \textit{Mathematische Annalen, }ejected Brouwer from the Editorial
Board.

\textit{Von Neumann.} Von Neumann contributed to foundational questions,
e.g. by formalising the (or a) construction of the natural numbers $\mathbb{N%
}$:

$0:=\{\emptyset \}$, $1:=\{0\}$, $2:=\{0,1\}$, $3:=\{0,1,2\}$, etc.: $%
n+1:=n\cup \{n\},$ see [Hal, \S 11] ([Neu1], [Neu2], 1928), and work on
amenable groups, with applications to the `Banach-Tarski paradox' (as above)
([TomW]; [Bin]).

The sets $x$ in Von Neumann's definition are ordered by $\in $ and are 
\textit{transitive}: if $z\in y\in x$, then $z\in x.$ Indeed the ordinals,
which form the class $On$ (not a set), are initially introduced as
transitive well-ordered \textit{structures} $\langle x,\in _{x}\rangle $
with $\in _{x}$ the restriction to $x$ of the membership relation. Once
ordinals $\alpha $ are established (this uses the axiom of regularity), the
cumulative hierarchy $V_{\alpha }$ may be introduced inductively so that $%
V_{\alpha +1}=\wp (V_{\alpha }),$ with $\wp $ the power set operation, and $%
V_{\lambda }=\bigcup \{V_{\alpha }:\alpha <\lambda \}$ for $\lambda $ a
limit ordinal. The class of sets is then $V=\bigcup \{V_{\alpha }:\alpha \in
On\},$ and each set $x$ has a well-defined \textit{rank}: the least $\alpha $
with $x\in V_{\alpha }.$

The formal language of set theory $LST$ builds formulas from a defined
sequence of free variables (e.g. $v_{0},v_{1},...),$ the atomic ones taking
the form $x\in y$ and $x=y,$ with $x$ and $y$ standing for free variables;
the syntactically more complex ones then arise from the usual logical
connectives and quantifiers ($\forall x$ and $\exists y$ -- creating bound
variables from the free variables $x,y).$ The idea is that the free and
bound variables are restricted to range only over the elements in the
universe of discourse (thus yielding a `first-order' language). This
language is a necessary ingredient of the axiomatic method, its first
purpose being to give meaning to the notion of `property' (so that e.g. $%
\{x\in y:\varphi (x)\}$ is recognized as a set when $\varphi $ is a formula
with one free variable $x$).

The language $LST$ is minimal as compared to the language of, say, group
theory, whose type (officially: `signature') involves more items (a
designated constant $1$, functions like $y\circ z$ , relations, etc). Each
such language is interpreted in a mathematical structure; for instance, at
its simplest a group structure has the form $\mathcal{G}:=\langle
G,1_{G},\circ _{G},\cdot ^{-1}\rangle $ and so lists its domain, designated
elements and operations. Below structures are assumed to be \textit{sets}
unless otherwise qualified; it is sometimes convenient (despite formal
complications) to allow a class as a domain, e.g. $\langle V,\in \rangle .$

The (metamathematical -- i.e. `external' to the discourse in the language)
semantic relation $\models $ of satisfaction/truth (below), due to Tarski
(see [Tar2], cf. [BelS, Ch. 3 \S 2]), is read as `models', or informally as
`thinks' (adopting a common enough anthropomorphic stance). A formula $%
\varphi $ of $LST$ with free variables $x,y,...,z$ may be interpreted in the
structure $\mathcal{M}:=\langle M,\in _{M}\rangle $ (with $\in _{M}$ now a
binary set relation on the set $M$) for a given assignment $a,b,...,c$ in $M$
for these free variables, and one writes%
\[
\mathcal{M}\models \varphi (x,y,...,z)[a,b,...,c],\text{ or by abbreviation }%
M\models \varphi \lbrack a,b,...,c] 
\]%
if the property holds; this requires an induction on the syntactic
complexity of the formula starting with the atomic formulas (for instance,
the atomic case $x\in y$ is interpreted under the assignment $a,b$ as
holding iff $a\in _{M}b$ ). Compare the reduction of complexity in the
forcing relation of \S 6 below.

This apparatus enables definition of `suitably qualified' forms of
`definability'; by contrast, unrestricted `definability' leads to such
difficulties as the `least ordinal that is not definable', so is to be
avoided (compare \S 3 below with Tarski's undefinability of truth). A simple
example is that of an element $w\in M$ being definable over $M$ from a
parameter $v\in M,$ in which case for some formula $\varphi (x,y)$ with two
free variables:%
\[
w\text{ is the unique }u\in M\text{ with }M\models \varphi (u,v). 
\]%
Thus G\"{o}del introduced the constructible hierarchy $L_{\alpha }$ by
analogy with $V_{\alpha }$: however, $L_{\alpha +1}$ comprises only sets
definable over $L_{\alpha }$ from a parameter in $L_{\alpha };$ here $%
L_{\lambda }=\bigcup \{L_{\alpha }:\alpha <\lambda \}$ for $\lambda $ a
limit ordinal, a matter we return to later, yielding the class $L=\bigcup
\{L_{\alpha }:\alpha \in On\}$.

Certain formulas, like $\varphi (x,y)$ above (which can be explicitly, and
so effectively, enumerated, as $\varphi _{m}$ say), may give rise via the
substitution of a parameter $v$ for $y$ to a family of not necessarily
unique elements $u\in M$ satisfying $\varphi (u,v).$ An appeal, in general,
to AC but in the `metamathematical' setting (i.e. the context of the
mathematics studying relations between the language and the structures),
selects a witness $w$ of the relation $\varphi (x,v)$ holding in $M$: the
function $v\mapsto w$ is called a \textit{Skolem function }(for\textit{\ }$M$
and $\varphi $); we will see a striking application presently -- for
background on this key notion see e.g. [Hod]. Evidently, a structure like $%
\mathcal{M}:=\langle L_{\alpha },\in _{L_{\alpha }}\rangle $ contains enough
well-orderings of its initial parts $L_{\beta }$ for $\beta <\alpha $
(induced by the enumeration $\varphi _{m}$ and well-ordering of the ordinal
parameters) that reference to AC\ here becomes unnecessary. (Incidentally,
this is why AC\ holds in the class structure $\langle L,\in \rangle $.)

We will refer to some other definability classes below in \S 6, so as an
introduction we mention two classical ones. The class $OD$ of ordinally
definable sets comprises those that are definable from ordinal parameters
over $\langle V_{\alpha },\in _{V_{\alpha }}\rangle $ for some $\alpha .$ An
element of a set in OD need not itself be in OD; the class HOD is the
smaller class of those elements $x$ whose transitive closure consists
entirely of sets in OD, so HOD is a transitive class; see [MyhS] for a
discussion.

In view of the finitary character of formulas, the L\"{o}%
wenheim-Skolem-Tarski theorem (see e.g. [Hod], or [BelS, Ch. 4.3]), as
applied to the language of set theory $LST$, asserts that if a set $\Sigma $
of sentences is modelled in a structure $\mathcal{M}$, then there exist
structures $\mathcal{N}$ of any infinite cardinality satisfying $\Sigma ,$
including \textit{countable} ones. The latter ones are generated by
induction by iterative application of all the Skolem functions; so this
needs only the Axiom of Dependent Choices. A familiar example is the
countable subring with domain $\mathbb{Q}$ of the ordered ring structure $%
\langle \mathbb{R},0,1,+,\times ,<\rangle $. Passing to above-continuum
cardinalities yields models of non-standard analysis with infinitesimals and
infinite integers (see below); but here AC is needed to construct Skolem
functions with which to generate the much larger structure.

The axioms of set theory include a finite set and an \textit{axiom schema}
corresponding to the Axiom of Replacement (which asserts that the image of a
set under a functional relation $\varphi (x,y)$ expressed in $LST$ is again
a set). In order to model these axioms in structures like $\langle M,\in
_{M}\rangle $ with $M$ a set, it is necessary to restrict attention to the
use of a finite number of instances of the axiom schema -- causing no
practical loss of generality, since any amount of mathematical argument will
necessarily do just that (for instance, a deduction of an inconsistency).
Thus, assuming the consistency of the axioms of set theory, any finite
subset of the axioms has a model $\mathcal{M}$ (by the G\"{o}del-Henkin
Completeness Theorem; see e.g. [BelS, Th. 4.2]) and so also a countable
model $\mathcal{N}$. This is conventionally and systematically rephrased as
saying that the axioms of set theory have a countable model; compare [Kun2,
Ch. 7 \S 9].

By its very nature the countable model $\mathcal{N}$ will contain far fewer
bijections than exist in Cantor's world $V.$ If transitive, the domain $N$
of $\mathcal{N}$ will have an initial segment of the ordinals in $V;$
however, there will be countable ordinals which $\mathcal{N}$ `thinks' are
uncountable, owing to missing bijections. The rule to observe is that 
\textit{ordinals are absolute} whereas \textit{cardinality is relative}.
This is exploited in arranging the failure of the Continuum Hypothesis, CH,
by the model extension process of forcing (see below for details and
references). In the context of a transitive model of set theory $\mathcal{M}$
we will write e.g. $\omega _{1}^{M}$ for the ordinal which in $\mathcal{M}$
is its first uncountable. In the absence of a superscript the implied
context is $V.$

Provided the Regularity axiom is included, the structure $\mathcal{N}%
=\langle N,\in _{N}\rangle ,$ being then well-founded, is isomorphic to a
transitive structure; the isomorphism $\pi $ is given inductively by:%
\[
\pi (x):=\{\pi (y):y\in _{N}x\}, 
\]%
and is known as the \textit{Mostowski collapse}. Thus, for example, $\pi
(\emptyset ^{M})=\emptyset .$

\bigskip

\textbf{3. G\"{o}del, Tarski and their legacy}

The use of formal language brought greater clarity to the axiomatic method:
thus Skolem helpfully clarified one of Zermelo's axioms by replacing the
latter's use of the informal notion of `definite property' with a formal
rendering (i.e. by reference to formulas in a formal language). This was
soon to be followed by the discovery of the limitations of formal language:
the publication in 1931 of G\"{o}del's two incompleteness theorems, preceded
by the results of his 1930 thesis on the completeness of first-order logic
(that every universally valid sentence is provable -- [BelS, Th. 12.1.3])
and on compactness (a corollary). The latter was to bear fruit at the hands
of Tarski much later (1958 on). We note that the Compactness Theorem for
predicate calculus (that a set of sentences has a model iff each finite
subset has a model [BelS, Ch. 5 \S 4]), Tychonov's Theorem in topology and
AC\ are deeply connected; see [Jec1]. See also [BelS, Ch. 5 esp. \S 5] for
the status of variants and the connection with the ultraproducts of \S 5
below.

The two incompleteness theorems concerning any axiomatic system rich enough
to encompass arithmetic (firstly, the existence in the formal language of
the axioms of sentences that can be neither proved nor disproved, and
secondly, the impossibility of such a system to provide a proof for its own
consistency), rather than just wreck Hilbert's programme, produced untold
benefits to the richness of mathematics: the plurality of the possible
interpretations of a set of axioms (as in Skolem's non-standard arithmetic),
and the accompanying search for choosing the ways to reduce incompleteness,
on the one hand, and to test or justify any belief in consistency, on the
other: especially in the case of the axioms of set theory. See [Ste].

G\"{o}del's enduring insight was the embedding by arithmetic coding (hence
the need for the `rich enough' presence of arithmetic) of (aspects of) a
`metalanguage' -- the informal language of discourse needed to examine a
formal language as a mathematical entity -- back into the formal language,
specifically the concepts of proof and provability -- see below.

Addressing the incompleteness of set theory, G\"{o}del's second legacy
relates to `relative consistency': proof in 1938 (published in 1940) of the
consistency relative to ZF of both AC -- a matter of supreme importance,
given the Banach-Tarski paradox (dating back to 1924) --\ and of GCH. The
key idea in the proof was the introduction (see \S 2 above) of the
cumulative hierarchy $L_{\alpha }$ of constructible sets whose totality
comprising the class $L$ is an \textit{inner model} (i.e. a subuniverse of
the universe $V\ $of von Neumann, specifically a transitive class containing 
$On$). This was to be the foundation stone for the advances of the `next one
hundred years' in two ways. The first was to invite extensions of $L$ by
appropriate choice of sets outside $L.$ The second, more technical, derives
from Skolem's method (1912) of constructing countable sub-models, enshrined
in a \textit{condensation principle,} that if $M\ $is a countable `submodel'
of $L$ (more accurately an `elementary substructure'), then it is isomorphic
to a set $L_{\alpha }.$

Contemporaneously with G\"{o}del's earliest contributions, and blending and
intertwining with them, there occurs a `volcanic eruption' of ideas and
results from the fertile mind of Tarski: bursting forth in 1924 with the
Banach-Tarski paradox (mentioned above) and evidenced by the working
seminars of 1927-1929, laying the foundations of Tarski's remarkable legacy,
both that published in its time and that published later. This included work
on the definability or otherwise (definable if `external', not if
`internal') of the concept of truth, a result closely allied to G\"{o}del's
incompleteness result and of similar vintage. Suffice it to point to the
role of `elementary substructure' (term due to Tarski) in the condensation
principle above.

Deficiencies in Hilbert's approach to geometry (e.g., its tacit assumption
of set theory) led Tarski to re-examine the axiomatic basis of geometry. In
1930 Tarski was able to prove the decidability of `elementary geometry', via
a reduction to `elementary algebra' where he was able to generalize Sturm's
algorithm for counting zeros of polynomials --- see [Vau] for references and
[SolAH] for recent developments in this area.

\bigskip

\textbf{4. Ramsey, Erd\H{o}s and their legacy: infinite combinatorics;
partition calculus and large cardinals}

\textbf{4a. Ramsey and Erd\H{o}s }Pursuing a special case of Hilbert's 
\textit{Entscheidungsproblem} of 1928 -- proposing the task of finding an
effective algorithm to decide the validity of a formula in first-order logic
-- Ramsey was led to results in both finite and \textit{infinite
combinatorics} (obtained late that same year, and published in 1930, [Ram]),
the finite version of which yielded the desired algorithm for the special
(\textquotedblleft though common\textquotedblright ) universal type of
formula. In general no computable algorithm exists, as was shown by Church
(using G\"{o}del's coding) in 1935, and independently by Turing in 1936 (via
Turing machines). The \textit{Infinite Ramsey Theorem} (which acted as a
paradigm for its finite variants) asserts in its simplest form that if the
distinct unordered pairs (doubletons) of natural numbers are partitioned
into two (disjoint) classes, then there exists an infinite subset $\mathbb{M}%
\subseteq \mathbb{N}$ all doubletons from which fall in the same class; thus 
$\mathbb{M}$, which may be said to be a \textit{homogeneous} (monochromatic)
subset for the partition, is large -- see [Dra, Ch. 2.8.1, Ch. 7.2 which
both use DC]. (Homogeneity is a constantly recurring theme in what follows.)
Thus, as a corollary, a Cauchy sequence in $\mathbb{R}$ contains either an
increasing or a decreasing subsequence. The combinatorial result extends
from doubletons to (unordered) $n$-tuples (called by Ramsey `combinations')
and from \textit{dichotomous} partitions to ones allowing any finite number $%
k$ of partitioning classes. Further analogues and generalizations form the
discipline of \textit{partition calculus}, the founding fathers of which
were Paul Erd\H{o}s and Richard Rado: see [ErdR].

Given its origins, it is not altogether surprising that Ramsey's theorem and
its generalizations continue to play a key role in the logical foundations
of set theory.

\bigskip

\textbf{4b. Partitions from large cardinals }We are particularly concerned
below with the partition property that follows. As usual we regard any
ordinal (including any cardinal) as the set of its predecessors. The
partition property (partition relation) of concern is%
\[
\kappa \rightarrow (\alpha )_{2}^{<\omega }, 
\]%
by which is meant that if $[\kappa ]^{<\omega }$ (the finite subsets of $%
\kappa $) is partitioned into two classes, then there is a homogeneous
subset of $\kappa $ of order type $\alpha .$ (Ramsey's result as stated
above is recorded in this notation as $\omega \rightarrow (\omega )_{2}^{2},$
and its immediate generalization to $n$-tuples and $k$ classes as $\omega
\rightarrow (\omega )_{k}^{n}.$)

For any $\alpha \geq \omega $ the least cardinal $\kappa $ for which $\kappa
\rightarrow (\alpha )_{2}^{<\omega }$ holds, denoted $\kappa (\alpha ),$ is
called the $\alpha $-th \textit{Erd\H{o}s cardinal }(or \textit{partition
cardinal}); but do such cardinals exist? One may show in ZFC that$\ \kappa
(\alpha ),$ if it exists, is regular (below), and when $\alpha $ is a limit
ordinal, that $\kappa =\kappa (\alpha )$ is strongly inaccessible (below)
[Dra, Ch. 10], and so $V_{\kappa }$ is a model of ZFC, written $V_{\kappa
}\models $ZFC. Hence, by G\"{o}del's incompleteness theorem, we cannot
deduce its existence in ZFC.

Of particular importance are cardinals $\kappa $, in particular $\kappa
=\kappa (\omega _{1})$, for which $\kappa \rightarrow (\omega
_{1})_{2}^{<\omega }$ holds: see the next section. So if, as we do, we need
them, then we must add their existence to our axiom system. To gauge the
consistency strength of this assumption we refer to one of the earliest
notions of a `large cardinal': a \textit{measurable cardinal} $\kappa .$
Such a cardinal was defined by Ulam [Ula] in 1930 by the condition that it
supports a $\{0,1\}$-valued $\kappa $-additive (i.e. additive over families
of cardinality $\lambda $, for all $\lambda <\kappa )$ non-trivial measure
on the power set $\mathcal{\wp }(\kappa )$. This may be reformulated as
asserting the existence of a $\kappa $-complete \textit{ultrafilter} on $%
\kappa $ ([Car2], [ComN], [Jec2], [GarP]). It turns out that for $\kappa $
measurable, the stronger relation $\kappa \longrightarrow (\kappa
)_{2}^{<\omega }$ holds. The latter is taken as the defining property of a 
\textit{Ramsey cardinal}, through its similarity with $\omega \rightarrow
(\omega )_{2}^{2}.$

We stop to notice that the relation $\kappa \rightarrow (\kappa )_{2}^{2}$
(taken to be the definition of a \textit{weakly compact cardinal} [Dra,
Ch.10.2]) holds iff $\kappa $ is strongly inaccessible and $\kappa $ has the
tree property: every tree of cardinality $\kappa $ having less than
cardinality $\kappa $ nodes at each level has a path, i.e. a branch of full
length $\kappa .$ It is interesting that, as with the Cauchy sequences in $%
\mathbb{R}$ above, if $\kappa \rightarrow (\kappa )_{2}^{2},$ then every
linearly ordered set of cardinality $\kappa $ has a subset of cardinality $%
\kappa $ which is either well-ordered or reversely well-ordered by the
linear ordering.

\bigskip

\textbf{4c. Large cardinals continued}

The first notion of a large cardinal is motivated by the conceptual leap
from the finite to the infinite, as exemplified by the set of natural
numbers viewed as $\mathbb{N}$, or, better for this context, as the first
infinite ordinal $\omega .$ The arithmetic operations of summation and
multiplication/exponentiation (equivalently, the power set operation $\wp )$
applied to members of $\omega $ lead to members below $\omega .$

This observation can be copied by a direct reference to the two
corresponding operations that generate a union of a given family and the
power set of a given set, each operation being guaranteed by the
corresponding axiom. Thus a cardinal is said to be \textit{weakly
inaccessible} if it is a limit cardinal above $\omega $ which is \textit{%
regular }(a regular limit cardinal), meaning, firstly, that it is the 
\textit{limit}, i.e. supremum (union), of all the preceding ordinals, and,
secondly, that nonetheless it is not the union (supremum) of a smaller
family of ordinals. A cardinal $\kappa $ is \textit{strongly inaccessible},
or just (plain) \textit{inaccessible}, if it is a regular strong limit
cardinal, i.e. additionally $2^{\lambda }<\kappa $ for all $\lambda <\kappa
. $ (Here $2^{\lambda }$ is the cardinality of $\wp (\lambda )$.) Further
such notions (of hyper-inaccessibility), which we omit here, have been
introduced by reference to the idea of a `large limit' (limit over a large
set) of `large cardinals'. The axioms ZFC, assumed consistent, cannot imply
the existence of an inaccessible $\kappa $, as then $V_{\kappa }$, being a
model for ZFC, provides proof within ZFC\ of the consistency of ZFC, a
contradiction to G\"{o}del's incompleteness theorem.

A second source of largeness is motivated by the study of infinitary
languages, the idea being to overcome some of the limitations of first-order
languages. For example, in the language $\mathcal{L}_{\kappa \kappa }$ one
admits $\kappa $ many free variables and permits infinite
conjunctions/disjunctions of a family of formulas of cardinality below $%
\kappa .$ This leads to the desirability of these languages having a
compactness property analogous to G\"{o}del's compactness property of the
ordinary language $\mathcal{L}_{\omega \omega }$ (see above). Examples of
the failure of compactness abound; so it emerges that the desired $\kappa ,$
if it exists, needs to be large. Thus a cardinal $\kappa $ is called \textit{%
strongly compact }[Dra, Ch. 10.3] if the language $\mathcal{L}_{\kappa
\kappa }$ is $(\lambda ,\kappa )$-\textit{compact} for each $\lambda \geq
\kappa $, that is: for each $\lambda \geq \kappa $ and any set $\Sigma $ of
sentences in that language with $|\Sigma |\leq \lambda ,$ if each subset $%
\Sigma ^{\prime }$ with $|\Sigma ^{\prime }|<\kappa $ has a model, then $%
\Sigma $ has a model. (So the cardinality of $\Sigma \ $here is not
constrained.) The property may be characterized without reference to the
language more simply as saying that every $\kappa $-complete filter can be
extended to a $\kappa $-complete ultrafilter.

Analogously, a cardinal $\kappa $ is \textit{weakly compact }[Dra, Ch. 10.3]
if the language $\mathcal{L}_{\kappa \kappa }$ is $(\kappa ,\kappa )$%
-compact: if any set of sentences $\Sigma $ with $|\Sigma |\leq \kappa $
such that each of its subsets of cardinality $<\kappa $ has a model, then $%
\Sigma $ has a model.

A third, more promising, source is more in keeping with the first
(`operational') viewpoint. It is motivated by the `substructures' analysis
initiated in G\"{o}del's proof that GCH\ holds in the universe of
constructible sets. Attention focusses now on the properties that the
operation of elementary embedding could or should have. We recall that the
range of such an embedding is an elementary substructure. Suppose that $%
j:N\rightarrow M$ is an elementary embedding, where $N$ and $\ M$ are
transitive classes and $j$ is definable in $N$ by a formula of set theory
with parameters from $N$. Then $j$ must take ordinals to ordinals and $j$
must be strictly increasing. Also $j(\omega )=\omega $ and $j(\alpha )\geq
\alpha ,$ so there is a least $\delta $ with $j(\delta )>\delta .$ This is
the \textit{critical point} of $j.$ Then%
\[
\mathcal{U}:=\{X\subseteq \delta :\delta \in j(X)\} 
\]%
is a non-principal $\delta $-complete ultrafilter on $\delta ,$ i.e. $\delta 
$ is a measurable cardinal. In fact, the converse is also true -- see
[SolRK, Th. 1.2]. Interestingly, here a non-principal ultrafilter is defined
by membership of a single point, albeit via images.

The significance of this characterization lies in the `operations' the
function $j$ encodes which, on the one hand, pass the test of `elementarity'
and, on the other, introduce an upward jump at the critical point (roughly
speaking, an `inaccessibility from below by elementarity').

We mention some further canonical large-cardinal notions obtained from
variations on this elementary embedding theme; these will be useful not only
presently for the establishement of a \textit{reference scale} of
consistency strength, but also later in relation to the regularity
properties of subsets of $\mathbb{R}$ (such as Lebesgue measurability etc.,
considered in \S 7 and 10).

A cardinal $\kappa $ is \textit{supercompact} if it is $\lambda $\textit{%
-supercompact} for all $\lambda \geq \kappa $; here $\kappa $ is $\lambda $%
-supercompact if there is a (necessarily non-trivial) elementary embedding $%
j=j_{\lambda }:V\rightarrow M$ with $M$ a transitive class, such that $j$
has critical point $\kappa $, and $M^{\lambda }\subseteq M$, i.e. $M$ is
closed under arbitrary sequences of length $\lambda $. Under AC, w.l.o.g. $%
j(\kappa )>\lambda $.

For $\kappa $ a cardinal and $\lambda >\kappa $ an ordinal, $\kappa $ is
said to be $\lambda $\textit{-strong} if for some transitive inner model (\S %
3), $M$ say, there exists an elementary embedding $j_{\lambda }:V\rightarrow
M$ with critical point $\kappa ,$ $j_{\lambda }(\kappa )\geq \lambda ,$ and%
\[
V_{\lambda }\subseteq M. 
\]%
Furthermore, $\kappa $ is said to be a \textit{strong cardinal} if it is $%
\lambda $-strong for all ordinals $\lambda >\kappa $.

This notion may be relativized to subsets $S$ to yield the concept of $%
\lambda $-$S$\textit{-strong }by requiring in place of the inclusion above
only that%
\[
j(S)\cap V_{\lambda }=S\cap V_{\lambda }. 
\]%
(One says that $j$ preserves $S\ $up to $\lambda $.) This provides passage
to our last definition. The cardinal $\delta $ is a \textit{Woodin cardinal}
if $\delta $ is strongly inaccessible, and for each $S\subseteq V_{\delta }$
there exists a cardinal $\theta <\delta $ which is $\lambda $-$S$-strong for
every $\lambda <$ $\theta .$

The \textit{consistency strength} of various extensions of the standard
axioms ZFC,\ by the addition of further axioms, may then be compared
(perhaps even assessed on a well-ordered scale) by determining which
canonical large-cardinal hypothesis will suffice to create a model for the
proposed extension. Thus, for $\kappa $ supercompact, $V_{\kappa }\models
\exists \mu \lbrack $\textquotedblleft $\mu $ is strong\textquotedblright $%
], $ which places supercompact above strong. Likewise, for $\kappa $ strong, 
$V_{\kappa }\models \exists \mu \lbrack $\textquotedblleft $\mu $ is
measurable\textquotedblright $]$, placing measurability below strong. (And
below that is the existence of a Ramsey cardinal, recalling earlier
comments.) The consistency of Woodin cardinals is thus between strong and
supercompact: diagramatically,%
\[
\text{supercompact \TEXTsymbol{>} Woodin \TEXTsymbol{>} strong \TEXTsymbol{>}
measurable \TEXTsymbol{>} Ramsey.} 
\]

\bigskip

\textbf{5. Beyond the constructible hierarchy }$L$ -- \textbf{I}

We have mentioned the L\"{o}wenheim-Skolem-Tarski theorem. How else may one
construct structures that will contain a given one as an elementary
embedding? In topology one naturally reaches for powers and products (as
with Tychonov's theorem), and also their various substructures such as
function spaces. For example, Hewitt [Hew] in 1948 constructed hyper-real
fields by using a quotient operation on the space of continuous functions
via a maximal ideal; cf. [DalW2].

\textbf{5a. Expansions via ultrapowers and intimations of indiscernibles}

Jerzy \L o\'{s} [\L o\'{s}] in 1955, though foreshadowed by Skolem's
construction [Sko] of non-standard arithmetic in 1934, and even G\"{o}del
1930, introduced a natural algebraic way of constructing new structures. \L o%
\'{s} relied on the concept, introduced in 1937 by Cartan [Car1,2], of 
\textit{ultrafilter}: a maximal filter in the power set of $I$, say. (The
assumption of the existence of these -- see PI in \S 1-- is in general
weaker than AC.) For a family of structures $\langle \mathcal{A}_{i}:i\in
I\rangle ,$ all of identical type/signature, i.e. each having the same
distinguished operations and relations on its domain $A_{i}$ (and possibly
distinguished elements), one first defines the direct product as a structure
(again of the same type) with domain the set $\tprod\nolimits_{i\in I}A_{i}$
(the product's existence in general implicitly invoking AC, of course) by
defining the operations and relations pointwise; thus any distinguished
element $e$, say, if interpreted in $A_{i}$ as $e_{i}$, say, is interpreted
in the product by the function $e:i\mapsto e_{i}$. Next, for $\mathcal{U}$
an ultrafilter on $I$, define $\mathcal{U}$-equivalence: $f\sim g$ according
as $\{i\in I:f(i)=g(i)\}\in \mathcal{U}$, i.e. $f$ and $g$ are pointwise $%
\mathcal{U}$-almost equal. Then denote by $\tprod\nolimits_{i\in I}\mathcal{A%
}_{i}/\mathcal{U}$ the equivalence classes $[f]_{\mathcal{U}}$ and equip
these with the requisite operations and relations suitably interpreted as
relations that hold pointwise $\mathcal{U}$-almost always.

By induction from the construction of these `atomic' cases of relations, \L o%
\'{s}'s Theorem (\L T below) asserts satisfaction in the ultraproduct of
general properties/formulas $\varphi ,$ say for simplicity with one free
variable $v,$ via 
\[
\dprod\nolimits_{i\in I}\mathcal{A}_{i}/\mathcal{U}\models \varphi (v)[f]_{%
\mathcal{U}}\text{ iff }\{i\in I:\mathcal{A}_{i}\models \varphi (f(i))\}\in 
\mathcal{U}, 
\]%
for $\varphi $ any first-order formula (in the language needed to describe a
structure of that type -- `signature' above).

If the $\mathcal{A}_{i}=\mathcal{A}$ are all equal (with domain $A),$ then $%
\mathcal{A}$ embeds elementarily into the \textit{ultrapower} $\mathcal{A}%
^{I}/\mathcal{U}$, when $a\in A$ is identified with the constant map $%
f_{a}:i\mapsto a$.

Consider $\mathcal{A}:=\langle \mathbb{R},+,\cdot ,\leq ,0,1\rangle ,$ $I=%
\mathbb{N}$ and $\mathcal{U}$ an ultrafilter extending the filter of
co-finite subsets of $\mathbb{N}$ (again invoking, say, AC). Then, $\mathbb{R%
}$ embeds in $\mathbb{R}^{I}/\mathcal{U}$ , with any real number $a$
represented by the constant function $f_{a}:$ $n\mapsto a.$ Let us call the
function $\mathrm{id}(i):=i$ for $i\in I$ a \textit{dominating function}
since it plays an important role and \textit{dominates} any constant
function $f_{m}$ for $m\in \mathbb{N}$; indeed, $[f_{m}]_{\mathcal{U}}\leq
\lbrack \mathrm{id}]_{\mathcal{U}}$, since $\{n:m\leq n\}\in \mathcal{U}$,
and so $\mathrm{id}$ is an element following all of $\mathbb{N}$, and so
follows all of $\mathbb{R}$ in $\mathbb{R}^{I}/\mathcal{U}$. That is, $%
\mathrm{id}$ identifies an \textit{infinite number}; likewise $1/\mathrm{id}$
identifies a positive (non-zero) element that may be interpreted as an%
\textit{\ infinitesimal}. (This observation allowed Abraham Robinson
[Rob1,2] to develop a \textit{non-standard analysis} within which to
interpret and interrogate rigorously Leibniz's intuitive texts on
infinitesimals; see [Kei] for an undergraduate rigorous development of
calculus in this setting.)

The argument just given may be repeated with $\mathcal{A}:=\langle A,\in
_{A}\rangle $ for $A$ a transitive set and $\in _{A}$ the relation of
membership in $A.$ If $\mathcal{A}$ is a countable model of ZF, then,
provided $\mathcal{U}$ is countably complete (see e.g. [Kan, Prop. 5.3]), $%
\mathcal{A}^{I}/\mathcal{U}$ is well-founded under its `interpretation of
the membership relation', so will contain elements that form an interval of
ordinals following the ordinals in $A.$ However, there are no means within $%
\mathcal{A}$ itself of `seeing' the existence of this extra layer of
ordinals: speaking informally (but see below), they are `indiscernible'.
(Strictly speaking, $\mathcal{A}^{I}/\mathcal{U}$ needs to be replaced by an
isomorphic structure which is a transitive set, known as the \textit{%
Mostowski collapse}$,$ defined inductively by the collapsing function $\pi $:%
\[
\pi ([f]_{\mathcal{U}})=\{\pi ([g]_{\mathcal{U}}):[g]_{\mathcal{U}}\in _{%
\mathcal{U}}[f]_{\mathcal{U}}\} 
\]%
(cf. \S 2); then interpretations of ordinals collapse to actual ordinals.)

When $I=\kappa $ with $\kappa $ the least measurable cardinal and $\mathcal{U%
}$ the ($\kappa $-complete) corresponding ultrafilter, Dana Scott considered
the extension of $L$ to $L[\mathcal{U}]$ (the L\'{e}vy class of sets
`constructible relative to' $\mathcal{U}$ -- obtained by allowing
definability over the ordinals to refer also to $\mathcal{U}$ -- so a class
closed under the intersection with $\mathcal{U}$; see [Kan, Ch. 1 \S 3],
[Dra, 5.6.2]), and investigated the ultrapower $L[\mathcal{U}]^{I}/\mathcal{U%
}$ to conclude the non-existence of a measurable cardinal in $L.$ This is
easiest to understand through the lens of the theorem that existence of a
measurable cardinal contradicts $V=L$ [Sco1], [Dra, 6.2.10], [BelS, Ch. 14 
\S 6] (so there is no measurable cardinal in $L$). This is done again by
referring to the \textit{dominating function} $\mathrm{id}(i)=i,$ which vies
with $\kappa $ for the place of smallest measurable cardinal (in the
Mostowski collapse). A proper proof needs to avoid doubtful manipulations of 
$\mathcal{U}$-equivalence classes of subclasses of $L[\mathcal{U}]$. (To
achieve this, one represents any function $f$ by one of least rank $\mathcal{%
U}$-equivalent to it -- the `Scott trick'; under these circumstances
well-foundedness of the resulting model needs to be verified, using $\sigma $%
-additivity of $\mathcal{U}$.)

The gist of the proof is to recreate the following contradictions stemming
from \L T. As before, $\mathrm{id}(i):=i$ for $i\in I,$ and $f_{\lambda
}:i\mapsto \lambda $ is the constant function on $I$ embedding $\lambda $
into the ultrapower. By \L T, the map $\lambda \mapsto f_{\lambda }$ is
injective for $\lambda <\kappa $ (since $\{i:f_{\lambda }(i)=\lambda <\mu
=f_{\mu }(i)\}=\kappa \in \mathcal{U},$ for $\lambda <\mu <\kappa ).$ By \L %
T again, $f_{\kappa }$ is the smallest measurable cardinal in $\mathcal{A}$
(since $f_{\kappa }(i)=\kappa $ is such a cardinal, for all $i),$ hence $%
f_{\kappa }=\kappa $ (up to equivalence, really). Now $\mathrm{id}<f_{\kappa
}$ (since $\{i\in \kappa :i=\mathrm{id}(i)<f_{\kappa }(i)=\kappa \}=\kappa
\in \mathcal{U}$), so $\mathrm{id}<f_{\kappa }=\kappa .$ But, for each $%
\lambda <\kappa ,$ we have $f_{\lambda }<\mathrm{id}$ (since $\{i:f_{\lambda
}(i)<\mathrm{id}(i)\}=\{i\in I:\lambda <i\}=\kappa \backslash (\lambda
+1)\in \mathcal{U}$, as $\mathcal{U}$ is $\kappa $-complete). But $%
\{f_{\lambda }:\lambda <\kappa \}$ has cardinality $\kappa ,$ and so $\kappa
\leq \mathrm{id},$ contradicting the earlier deduction that $\mathrm{id}%
<\kappa .$

Actually, these observations just demonstrate that the embedding $j=j_{%
\mathcal{U}}$ obtained by composing $\lambda \mapsto f_{\lambda }$ with the
Mostowski collapse satisfies $j(\lambda )=\lambda $ for $\lambda <\kappa ,$
and $j(\kappa ),$ being the collapsed version of $[\mathrm{id}],$ lies
strictly above $\kappa ;$ thus the ordinal $\kappa $ is the \textit{critical
point} of $j.$

This argument was further investigated by Haim Gaifman, from the point of
view of iterating the ultrapower construction, and perfected by Kunen [Kun1].

\bigskip

\textbf{5b. Ehrenfeucht-Mostowski models: expansion via indiscernibles}

At about the same time as \L o\'{s} introduced ultraproducts into
model-theory, Ehrenfeucht and Mostowski [EhrM] in 1956 introduced a
construction that expands a structure $\mathcal{A}$ by importing a linearly
ordered set of elements in such a way that, speaking anthropomorphically, $%
\mathcal{A}$ is incapable of distinguishing between these imports and a
certain infinite subset of its own domain. Less than a decade later, first
Morley in 1962 (see e.g. [Mor]) and then Silver in his thesis in 1966 (see
[Sil]) put these features to decisive use, by enabling the imported elements
to generate various kinds of information about $\mathcal{A}$ consistent with
that generated by $\mathcal{A}$ on its own.

The original construction provided an elementary embedding of any infinite
structure $\mathcal{A}$ into another `larger' one -- larger in possessing
many non-trivial automorphisms, securing in particular a non-trivial
elementary embedding. A (copy of a) linearly ordered set $X$ is adjoined to $%
A$ of elements $x$ which are to be `indiscernible' from the viewpoint of $%
\mathcal{A}$ (except only in name -- as the formal language must adjoin
formal names $c_{x}$ to speak about them) in the sense that:%
\[
(\mathcal{A},(c_{x})_{x\in X})\models \varphi
(x_{1},...,x_{n})\Leftrightarrow \varphi (x_{1}^{\prime },....,x_{n}^{\prime
}), 
\]%
for all formulas $\varphi $ having $n$ free variables, for all $n,$ and all $%
x_{1}<...<x_{n}<x_{1}^{\prime }<...<x_{n}^{\prime }$ in $X.$

That this is possible in general relies on the Compactness Theorem (and so
on AC): the idea here being that if one takes the sentences true in $%
\mathcal{A}$ together with sentences $\varphi
(c_{x_{1}},...,c_{x_{n}})\Leftrightarrow \varphi (c_{x_{1}^{\prime
}},...,c_{x_{n}^{\prime }})$ (also the inequalities $c_{x}\neq c_{y}$), then
one may satisfy a finite set $F$ of these by interpreting the finite number $%
m$ of $c_{x}$s in play in $F$, $c_{x_{1}},...,c_{x_{m}}$ say, with suitably
chosen elements of $A,$ as follows. To effect the choice, partition all $m$%
-tuples of $A$ according as to whether or not $\mathcal{A}$ can distinguish
between them on the basis of the properties defined by the finite number of
formulas $\varphi (v_{1},...,v_{m})$ obtained from the $\varphi $ in $F.$
(That is: the free variables $v_{i}$ replace the constants $c_{x_{i}}$.)
Then an infinite homogenous set for this partition yields a model for $F.$

In particular, for limit ordinal $\delta ,$ the structure $\mathcal{A}%
=\langle L_{\delta },\in \rangle $ (by abuse of notation $\in $ here and
below denotes membership $\in $ restricted to $L_{\delta }$) can be expanded
to a structure with a sequence of indiscernibles whose formal language names
are $c_{n}$. Call that $\mathcal{A}_{0}.$ (Here AC may be avoided, as $%
L_{\delta }$ is well-ordered.) In turn, for any ordinal $\alpha ,$ that
expanded structure $\mathcal{A}_{0}$ may be further extended to a structure $%
\mathcal{M}_{\alpha }(\mathcal{A})$ with a set of indiscernibles $X$ of
order type $\alpha $ and with the following additional property: for any
formula in the language of $\langle L_{\delta },\in \rangle ,$ $\varphi
(v_{1},...,v_{n})$ say, 
\[
(\mathcal{A},(c_{n})_{n\in \omega })\models \varphi (c_{1},...,c_{n})\text{
iff }\mathcal{M}_{\alpha }\models \varphi (\mathbf{x})\text{ for some }%
\mathbf{x}=(x_{1},...,x_{n})\in X^{n}. 
\]%
So, in particular, the indiscernibles $X$ can generate all the true
sentences about $\mathcal{A}$. But are the structures $\mathcal{M}_{\alpha }(%
\mathcal{A})$ well-founded for \textit{all} $\alpha $? That depends on
whether the structures $\mathcal{M}_{\alpha }(\mathcal{A})$ for just $\alpha
<\omega _{1}$ are all well-founded (the reduction here is possible, since
any descending sequence occurring in the models with larger $\alpha $ can be
captured by a countable submodel). This will be so when $\mathcal{A}=\langle
L_{\kappa },\in \rangle $ and $\kappa $ satisfies the partition relation%
\[
\kappa \rightarrow (\omega _{1})_{2}^{<\omega }. 
\]%
(With $\alpha <\omega _{1}$ as above, the argument is similar to but easier
than that in the Ehrenfeucht-Mostowski result. Appealing to the partition
relation above in place of Ramsey's theorem, partition $(\xi _{1},...,\xi
_{n})\in \lbrack \kappa ]^{<\omega }$ dichotomously according as to whether $%
\mathcal{M}_{\alpha }\models \varphi (\xi _{1},...,\xi _{n})$ holds or not;
extract an $\omega _{1}$ homogeneous subset of $\kappa $ and use its first $%
\alpha $ members as the required indiscernibles. Their Skolem hull in $%
L_{\kappa }$, a well-founded set, is isomorphic to $\mathcal{M}_{\alpha }(%
\mathcal{A})$.)

A first corollary (by appeal to indiscernibility, use of only the first $%
\omega $ indiscernibles, and then the countability of the formal language):
only a countable number of subsets of $\omega $ are constructible in $L,$
even though from the viewpoint of $L$ there are uncountably many of them in $%
L$; but then, an embellishment of the analysis yields that $\omega _{1}^{L},$
the ordinal intepreted by $L$ as the first uncountable, is also countable.

Silver deduced deeper results about $L$ along these lines. Some of these
were then bettered by Kunen [Kun1], who devised a way for iterating the
ultrapower construction of a structure $\mathcal{M}$ in a setting where the
ultrafilter $\mathcal{U}$ need not be a member of $\mathcal{M}$. A most
remarkable contribution from Silver was the introduction of the set now
called $0^{\#}$ (\textit{zero-sharp}) following Solovay (originally
designated a `remarkable' set); this is the set of G\"{o}del codes $\lceil
\varphi \rceil $ for all the true sentences $\varphi $ about $L$ generated
by the $\omega $-sequence of indiscernibles $\{\omega _{1},\omega
_{2},...,\omega _{n},...\},$ namely: 
\[
0^{\#}:=\{\lceil \varphi \rceil :L\models \varphi (x_{1},...,x_{n})\text{
for }(x_{1},...,x_{n})\in \{\omega _{1},\omega _{2},...,\omega _{n},...\}\}. 
\]%
(The notation tacitly assumes that $n=n(\varphi )$ is the number of free
variables in $\varphi .)$ This set's very existence of course depends on
suitable large-cardinal assumptions, such as $\kappa \rightarrow (\omega
_{1})_{2}^{<\omega }$ holding for some $\kappa .$ The `existence of $0^{\#}$%
' can be used as a large-cardinal assumption in its own right, lying below
the existence of the Erd\H{o}s cardinal. Indeed, in \S 7 we discuss the
classical theory of analytic sets and thereafter the determinacy of infinite
positional games with a target set $T,$ say; the assumption that sets with
co-analytic target set are determined ($\mathbf{\Pi }_{1}^{1}$-determinacy)
implies that 0$^{\#}$ exists, a result due to Harrington [Har].

We return to the indiscernibles for the structures $\mathcal{A}=\langle
L_{\delta },\in \rangle $, \textit{assuming the partition relation} just
mentioned, which had been studied initially by Gaifman and by Rowbottom.
Silver's great contribution was to describe the structure, indeed the `very
good behaviour' (below), of a (proper) class $X$ of ordinal indiscernibles:
closed (under limits -- i.e. under suprema), unbounded in any cardinal $%
\lambda $ (with $X\cap \lambda $ of cardinality $\lambda )$; with $L_{\alpha
}\prec L_{\beta }$ for $\alpha <\beta $ both in $X$ (indeed, stretching the
notation to class structures, with $L_{\alpha }\prec L$); having the
property that every set in $L$ is definable from parameters in $X.\ $Among
the significant consequeness is the, already mentioned, countability of
those sets in $L\ $that are definable over $L$ without any parameters
(implying immediately that $V\neq L),$ and more importantly the definability
of truth in $L.$ For details see e.g. [Dra, Th. 4.8]. We stress these
results are subject to the partition assumption.

The point (above) about good behaviour concerns particularly the `closed
unbounded' nature of $X$ above. Sets of ordinals with this property should
be regarded as `large', since they enable the very important `stationary
sets' of the next section to be thought of as non-negligible. The two
concepts play a leading role in combinatorial principles (holding in $L)$
isolated by Jensen (see e.g. [Dev1]) from the fine structure of $L$. These
include Jensen's $\Diamond $ (diamond), used in constructing a `Suslin
continuum' as a counterexample to Suslin's hypothesis (see below); $\square $
(square); derived ones like $\clubsuit $ (club), introduced by Ostaszewski
[Ost1] (in `counterexample' constructions for general topology); and
generalizations $\clubsuit _{\text{NS}}$ studied by Woodin [Woo1, Ch. 8].
Compare the use of NT (for No Trump) in [BinO1, 4].

\bigskip

\textbf{6. Beyond the constructible hierarchy }$L$ -- \textbf{II}

\textbf{6a. Forcing and generic extensions}

The undisputed game-changer for set theory was Cohen's `method of forcing',
devised as a means of importing into a countable structure $\mathcal{M}%
=\langle M,\in _{M}\rangle $ additional sets from $V\backslash M$ ($V$
contains the reals; $M,$ being countable, does not), without disturbing the
fact that $\mathcal{M}$ may be a model of ZF. Speaking anthropomorphically,
the imported set may have the intention of introducing new information --
say, the existence of a transfinite sequence of real numbers viewed by $%
\mathcal{M}$ as an $\omega _{2}^{M}$ sequence (reference here to the
interpretation in $\mathcal{M}$ of the second uncountable cardinal), albeit
viewed by $V$ as a countable sequence -- without nevertheless encoding such
catastrophic information as that $M$ itself is countable. Cohen described
his method [Coh3] as ultimately analogous to the construction of a field
extension: introduce a name for the algebraically absent element, and then
describe its properties via polynomials in that element. In truth the
extension method shares a family resemblance with non-constructive existence
proofs, either via the Baire category method (the desired item has generic
features), or the Erd\H{o}s probabilistic method (measure-theoretic: the
desired item has `random' features). Indeed the two canonical instances of
forcing to adjoin real numbers, Cohen's and Solovay's, are categorical
(Cohen reals) or measure-theoretic (`random reals', or -- perhaps better --
`Solovay reals'). Indeed, following an idea of Ryll-Nardzewski and of
Takeuti, Mostowski [Most] shows how to guide the selection of an imported
set by reference to the points of a Baire topological space (one in which
Baire's theorem holds); avoiding a specified meagre set ensures that the
extension of $\mathcal{M}$ will be a model of ZF. The two canonical cases
then correspond to two topological spaces. For an alternative unification
see [Kun3].

One views the forcing method as acting `over' a structure $\mathcal{M}$ by
providing a set $P$ in $M$ of partial descriptions of a generic object $G$
yet to be determined. $P$ is thus rendered as a partially ordered set, and
under its ordering relation $q\leq p$ is understood as saying that $q$
contains more information about the object to be constructed than does $p.$
There is a syntactic relation $p\Vdash \varphi $ for $p\in P$ and $\varphi $
a sentence, read as `$p$ forces $\varphi $', which may be `explained' by an
induction reminiscent of the Tarski inductive definition of truth ($\models
, $ in \S 2), but with significant differences (below).

Before embarking on the details, it is helpful to use an analogy with
probability or statistical inference. Indeed $p\in P$ is usually called a
`condition'; forcing is inspired by the language of `conditioning'; its
inferences are concerned with information about $G$ given the information in 
$p.$ Thus the forcing relation must allow for further information which may
become available `later', so to speak.

As a first pass, here is a brief glimpse of the character of the forcing
relation: as this is a syntactical relation, we refer to a language whose
terms are built from functions from $P$ to $M,$ and so we have (see [Kun2,
Cor. 3.7]): 
\begin{eqnarray*}
p &\Vdash &\varphi \wedge \psi \text{ iff }p\Vdash \varphi \text{ and }%
p\Vdash \psi , \\
p &\Vdash &\lnot \varphi \text{ iff not}(\exists q\leq p)\text{ [}q\Vdash
\varphi \text{],} \\
p &\Vdash &(\exists v)[v\in \sigma \wedge \varphi (v)]\text{ then }(\exists
q\leq p)(\exists x\in \text{dom}(\sigma ))[\text{ }q\Vdash \varphi (x)].
\end{eqnarray*}

A clearer picture will emerge shortly.

Whilst a variant of the forcing relation above was Cohen's starting point,
this is now a derived concept, the usual starting point being a set $G$ that
is a filter on $P$ with the property that whenever $D\ $is a \textit{dense}
subset of $P$ (i.e. for each $p$ there is $q\leq p$ with $q\in D)$ and $D\in
M,$ then%
\[
G\cap D\neq \emptyset . 
\]%
Then $G$ is said to be $P$-\textit{generic }over $\mathcal{M}$, or just 
\textit{generic }over $\mathcal{M}$, when $P$ is understood.

For $M$ countable, the dense subsets of $P$ lying in $M$ may be enumerated
as a sequence $D_{n},$ and we may choose $p_{n}\in P$ starting with an
arbitrary $p_{0}\in D_{0}$ and inductively $p_{n+1}\leq p_{n}$ with 
\[
p_{n+1}\in D_{n+1}. 
\]%
The choice is possible precisely because $D_{n+1}$ is dense. Then $%
G:=\{q:(\exists n)$ $q\leq p_{n}\}$ meets each $D_{n},$ and so is generic
over $\mathcal{M}$. This construction is sometimes called the \textit{Cohen
diagonalization argument,} since, in particular, $G$ decides every sentence $%
\varphi $. Indeed, the following set is dense:%
\[
D_{\varphi }:=\{p:p\Vdash \lnot \varphi \text{ or }p\Vdash \varphi \} 
\]%
(as $p\notin D_{\varphi }$ implies not($p\Vdash \lnot \varphi )$ and so $%
(\exists q\leq p)$ [$q\Vdash \varphi $]).

The idea is that the dense sets provide a structured way of hinting at the
properties of $G$ and about the various ways that $G$ might be selected, but
conditional on some given state of knowledge $p.$ The sequence $p_{n}$ above
runs through all possible dense sets in an arbitrary order, and brings into
existence a particular realization of $G.$

Before $G$ is created there are only names for $G$ and for all the possible
objects in the intended extension, given simply by the \textit{functions }in 
$M^{P}.$ (This corresponds to the use of \textit{polynomials} in field
extension.) But, once a generic $G$ is given, one may proceed inductively to
give an interpretation $\tau ^{G}$ to the `names' $\tau \in M^{P}$ of
objects, inductively so that%
\[
\tau ^{G}:=\{\sigma ^{G}:(\exists p\in G)[\sigma =\tau (p)]\} 
\]%
(mirroring the Mostowski collapse above), and so construct the extension $%
\mathcal{M}[G]$ as the set of $G$-interpretations. One then defines forcing
(relative to $P$ and $\mathcal{M}$) by:%
\[
p\Vdash \varphi \text{ iff (}\forall G\text{ generic over }\mathcal{M}\text{%
)[}p\in G\rightarrow \mathcal{M}[G]\models \varphi \text{].} 
\]%
This should clarify the three properties of the forcing relation introduced
earlier.

It emerges that if $\mathcal{M}\models $ ZFC, then $\mathcal{M}[G]\models $
ZFC. Furthermore, if $P$ satisfies the so-called \textit{countable chain
condition }(\textit{`ccc'})\textit{\ }(which actually calls for antichains
of $P$ in $M\ $to be countable), then all ordinals that are cardinals from
the viewpoint of $\mathcal{M}$ continue to be cardinals from the viewpoint
of $\mathcal{M}[G],$ and their cofinalities [Jec2] remain the same.

To secure the failure of CH, Cohen used as his conditions finite sets $p$
with elements of the form:%
\[
\langle n,\alpha ,i\rangle \text{ for }n\in \omega ,\text{ }\alpha <\omega
_{2},\text{ }i\in \{0,1\}, 
\]%
which act as coded messages about objects, named as $c_{\alpha },$ to be
imported from outside $M$ asserting that $n\notin c_{\alpha }$ if $i=0$ and $%
n\in c_{\alpha }$ if $i=1.$ As with the `dog that did not bark', that which $%
p$ will never say allows us to infer that $c_{\alpha }$ will be a subset of $%
\omega :$ this is \textit{forced} to be the case, since no extension of the
coded message $p$ can say otherwise. Thus $p$ `hints at information' by the
absence of information.

Formally, the corresponding $P,$ called $Add(\omega ,\omega _{2})$ since it
adds $\omega _{2}$ many subsets of $\omega $, may be defined in $M$ to
comprise `partial functions' $p$ with finite domain contained in $\omega
\times \omega _{2}^{M}$ and range in $\{0,1\}$, and with the ordering of
`increasing informativeness' that $q\leq p$ if $p\subseteq q,$ that is, $q$
contains at least all of the information in $p.$ The filter $G$ in $P\ $has
the property that $\tbigcup G=\{\langle n,\alpha ,i(n,\alpha )\rangle :$ $%
n\in \omega ,\alpha \in M\cap \omega _{2}^{M}\}$ for some $i:(n,\alpha
)\mapsto \{0,1\}.$ Indeed, for $n,\alpha $ as above, each of the sets%
\[
D_{n,\alpha }:=\{p:\langle n,\alpha ,i\rangle \in p\text{ for some }i\in
\{0,1\}\}
\]%
is dense, as may be readily checked. (Hint: Given $p\notin D_{n,\alpha }$
choose $q$ to contain both $p$ and $\langle n,\alpha ,1\rangle .)$ So $G\ $%
must meet $D_{\alpha ,n}$ for each $\alpha \in M$ (as $\omega \subseteq M,$
since $\mathcal{M}\models $ ZFC$).$ For $\alpha \in M\cap \omega _{2}^{M},$
put%
\[
G_{\alpha }:=\{n\in \omega :\langle n,\alpha ,i(n,\alpha )\rangle \in
G,i(\alpha ,n)=1\}\subseteq \omega .
\]%
Moreover, for distinct $\alpha ,\beta <\omega _{2}$, put%
\[
\Delta _{\alpha ,\beta }:=\{p:\langle n,\alpha ,i\rangle ,\langle n,\beta
,1-i\rangle \in p\text{ for some }n\in \omega \text{ and some }i\in
\{0,1\}\},
\]%
which is dense. (Given $p\notin \Delta _{\alpha ,\beta }$ choose $q$ to
contain both $p$ and $\langle m,\alpha ,1\rangle ,\langle m,\beta ,0\rangle $
for large enough $m.)$ So for distinct $\alpha ,\beta \in M\cap \omega
_{2}^{M},$ $G$ contains $\langle n,\alpha ,i\rangle ,\langle n,\beta
,1-i\rangle $ for some $n$ and $i,$ with $i=1,$ say (w.l.o.g.). Then $n\in
G_{\alpha }\backslash G_{\beta }$. Thus in $\mathcal{M}[G]$ there are $%
\omega _{2}^{M}$ distinct subsets of $\omega ,$ and so from the viewpoint of 
$\mathcal{M}[G]$ the continuum is at least $\omega _{2}$ (since $\omega
_{2}^{M}$ is still the interpretation of $\omega _{2}$ in $\mathcal{M}[G]$
by the ccc, which is satisfied by $P$ here).

We have just given an example of importing a set in order to increase the
cardinality of the continuum. (Note that this construction may be repeated
with $\omega _{2}^{M}$ replaced by $\omega _{\tau }^{M}$ for $\tau $ with
any cofinality other than $\omega $, that being the only restriction on the
cofinality of the continuum.)

An important ingredient in Solovay's result [Sol3] on LM (in constructing a
model of ZF+DC in which all sets of reals are Lebesgue measurable -- cf.
[Kan, Ch.13 \S 11] -- to which we refer in \S 10.2) uses a further partial
order $P^{\kappa }=Coll(\omega \times \kappa ,\kappa )$, introduced by L\'{e}%
vy, whose function is to alter/collapse a (strongly) inaccessible cardinal $%
\kappa $ so that in the extension $\mathcal{N}=\mathcal{M}[G^{\kappa }]$ ($%
G^{\kappa }$ being $P^{\kappa }$-generic over $\mathcal{M}$) it is the
ordinal $\kappa $ that appears as the first uncountable cardinal $\omega
_{1}^{N}$. Consequently the ordinals below $\kappa $ are made to be
countable by the importation of appropriate enumerations. Interest focuses
on the substructure $\mathcal{N}_{1}$ with domain the sets that are
hereditarily definable over $\mathcal{N}$, from a parameter in $\mathcal{%
N\cap }On^{\omega }$ (i.e. from a sequence of ordinals in $\mathcal{N}$),
much as defined earlier. $\mathcal{N}_{1}$ satisfies the axioms ZF (see
[MyhS]), and, significantly here, shares the same sequences of ordinals, in
particular \textit{the same reals}. (Here the reals are identified via
binary expansion with characteristic functions of subsets of $\omega $.)

The L\'{e}vy conditions (elements of $P^{\kappa })$ this time are partial
functions with finite domain $\omega \times \kappa $ and range in $\kappa .$
Since there are no bounds placed on the range values of the partial function
in this $P$, it follows that for $\alpha <\kappa $ the functions 
\[
G_{\alpha }:=\{\langle n,\lambda (\alpha ,n)\rangle :\langle \alpha
,n,\lambda (\alpha ,n)\rangle \in G^{\kappa }\} 
\]%
will collectively witness (by enumeration) that each $\lambda <\kappa $ is
countable. This ensures that $\kappa $ \textquotedblleft viewed
from\textquotedblright\ $\mathcal{M}[G^{\kappa }]$ is $\omega _{1}.$ Solovay
's purpose is to turn any transfinite sequence of ordinals below an
inaccessible $\kappa $ into an $\omega $-sequence. This helps him turn an
arbitrary set of reals $A$ that lies in $\mathcal{N}_{1},$ initially
definable in $\mathcal{N}$ via ordinal parameters, into one that is
definable via a real $a$. (This also carries the advantage that, since $%
\kappa $ retains its inaccessibility in the extension $\mathcal{M}[a]$, one
may w.l.o.g. argue as though $\mathcal{M}[a]$ is $\mathcal{M}.)$ As both $%
\mathcal{N}$ and $\mathcal{N}_{1}$ have the same reals, they also have the
same Borel sets and the same \textit{null} $\mathcal{G}_{\delta }$-sets.

Solovay's surprising innovation was to force over $\mathcal{M}[a]$ using the
non-null Borel sets $\mathcal{B}_{+}$ ordered by inclusion (smaller sets
yielding more information as to location). The key idea here is to introduce
the notion of a \textit{random real}, namely a real that cannot be covered
by any null $\mathcal{G}_{\delta }$-set coded canonically by a real $c$ of
the model $\mathcal{M}[a]$. (Solovay thought of these as `random'; we have
already mentioned that Cohen reals are categorical, while random (`Solovay')
reals are measure-theoretic; the term generic was already in use, so
unavailable. Compare our earlier use of the language of probability and
statistical inference above. One might also mention the term pseudo-random
number in computer simulation.) But, $\mathcal{M}[a]$ being countable, there
are only countably many such codes, so in $V$ the set of non-random reals is
null. For a set $A\subseteq \omega ^{\omega }$ that is definable from an $%
\omega $-sequence of ordinals (i.e., by a sequence from $On^{\omega }$),
suppose that with $a$ as above, for some formula $\varphi _{A}$ say, $A=\{x:%
\mathcal{M}[a][x]\models \varphi _{A}[a,x]\}$. Now one may choose a formula $%
\psi _{A}$ such that, for $G_{+}$ a $\mathcal{B}_{+}$-generic filter and any 
$x\in \mathcal{M}[G_{+}]\cap On^{\omega },$ 
\[
\mathcal{M}[G_{+}]\models \varphi _{A}[a,x]\text{ iff }\mathcal{M}[x]\models
\psi _{A}[a,x]. 
\]%
In $\mathcal{B}_{+}$ choose a maximal (necessarily countable, by positivity
of measure here) antichain of Borel sets $\mathcal{C}$ whose elements
`decide' the formula $\psi _{A}[\check{a},\dot{r}]$ (i.e. force the formula
or its negation), where $\check{a}$ is a name for the set $a$ given above,
and $\dot{r}$ is a name for a random real (cf. the use of $\dot{q}$ in \S 6b
below). Then for $x$ random:%
\[
x\in A\Longleftrightarrow x\in \tbigcup \{F_{c}:F_{c}^{\mathcal{M}[G]}\in 
\mathcal{C}\text{ and }F_{c}^{\mathcal{M}[G]}\Vdash \varphi _{A}[\check{a},%
\dot{r}]\}. 
\]%
Here $F_{c}$ is a non-null closed set canonically coded by $c.$ So modulo
the null set of non-random reals, $A$ is an $\mathcal{F}_{\sigma }.$

\textbf{6b. Forcing Axioms}

Solovay's argument makes heavy use in various ways of `two-step extensions'
like $\mathcal{M[}G][H]$ with $G$ an $\mathcal{M}$-generic filter and $H$ an 
$\mathcal{M[}G]$-generic filter. By implication, $G$ is associated with a
partial order $P$ in $\mathcal{M}$ and $H$ with a partial order $Q$ in $%
\mathcal{M[}G].$ This can be turned into a one-step extension $\mathcal{M[}%
K] $, but in a perspicuous way (more general than cartesian products), so
that a generic extension of a generic extension is again a generic
extension. Since the model $\mathcal{M[}G]$ is created by interpreting
`names' (using $G $ as in $\tau ^{G}$ above), the partial order for the
equivalent single step needs to be built out of $P$ and out of a name $\dot{Q%
}$ for $Q,$ and must refer to pairs $(p,\dot{q})$ with $p\in P$ and $\dot{q}$
a name for something that is $P$-forced to lie in $\dot{Q};$ likewise, the
order on the resulting composition of the two partial orders, denoted $%
P_{1}\ast \dot{Q}$, must make use of how the $P$-conditions $P$-force the
extension property $\dot{q}\leq \dot{q}^{\prime }$ between names for
elements of $\mathcal{M[}G][H].$ Thus a kind of syntactical analysis in $%
\mathcal{M}$ underlies this `iterated forcing'. More generally, any ordinal $%
\alpha $ of $\mathcal{M}$ can provide the basis for $\alpha $-step
iterations, and, as with the topologies on products so too here, various
kinds of $\alpha $-iterations may be constructed by appropriate constraints
on the supports (e.g. finite or countable). We omit the details, except to
mention that it was by use of such an iteration that Solovay and Tennenbaum
[SolT] showed that it is consistent that no Suslin continuum exists (so
otherwise than in $L,$ where such exists); this led to the more general
observation, proved by Martin and Solovay: the consistency of Martin's
Axiom, MA ([MarS], cf. [Fre1]), namely the statement that for all cardinals $%
\kappa $ below the continuum ($\kappa <\mathfrak{c)}$ the following holds:

\bigskip

$\mathrm{MA}(\kappa ):$ \textit{for every partial order }$P$\textit{\
satisfying the countable chain condition (ccc), and any family} $\mathcal{F}$
\textit{with} $|\mathcal{F}|\leq \kappa $ \textit{of dense subsets of }$P,$%
\textit{\ there is a filter }$G$\textit{\ in }$P$\textit{\ which meets each }%
$D\in \mathcal{F}$.

\bigskip

The reader will notice the similarity between the property of $G$ here and
that of a filter $P$-generic over $\mathcal{M}$; indeed Martin (and
independently Rowbottom) proposed this axiom as a combinatorial principle
that is forcing-free -- so, in particular, with the potential for immediate
applicability without expertise in logic. That potential was so quickly
realized both in theorem-proving and counterexample-manufacture -- look no
further than [Fre1] -- that it became the `tool of first choice' when
abstaining from CH whilst harbouring CH-like intuitions, because, like
Zorn's Lemma, it encapsulates a `construction without (transfinite)
induction', replacing the latter with a side-condition swept away into $%
\mathcal{F}$, the family of dense sets. Of course, the `implied' induction
was performed, off-line so to speak, in the Martin-Solovay paper [MarS],
aptly titled `Internal Cohen extensions', reflecting the view that MA
asserts that the universe of sets is closed under a large class of generic
extensions.

In regard to MA's huge significance as an alternative to the continuum
hypothesis: we cite after Martin and Solovay [MarS] the statistic that at
least 71 of 82 consequences of CH, as given in Sierpi\'{n}ski's monograph
[Sie], are decided by MA or $[$MA \& $2^{\aleph _{0}}>\aleph _{1}].$ Amongst
these are that MA implies:

\noindent (1) $2^{\aleph _{0}}$ is not a real-valued measurable cardinal;%
\newline
\noindent (2) the union of less than $2^{\aleph _{0}}$ (Lebesgue) null
/meagre sets of reals is null/meagre;\newline
\noindent (3) Lebesgue measure is $2^{\aleph _{0}}$-additive;\newline
and that $[$MA \& $2^{\aleph _{0}}>\aleph _{1}]$ implies: \newline
\noindent (1) Suslin's hypothesis, that every complete, dense, linear order
without first and last elements in which every family of disjoint intervals
is at most countable (the Suslin condition) is order-isomorphic to $\mathbb{R%
}$;\newline
\noindent (2) every $\mathbf{\Sigma }_{2}^{1}$ set of reals (for the $%
\mathbf{\Sigma }$ and $\mathbf{\Pi }$ notation of the projective hierarchy
see \S 9) is Lebesgue measurable and has the Baire property;\newline
\noindent (3) every set of reals of cardinality $\aleph _{1}$ is $\mathbf{%
\Pi }_{1}^{1}$ (co-analytic) iff every $\aleph _{1}$ union of Borel sets is $%
\mathbf{\Sigma }_{2}^{1}$.

\bigskip

It is worth remarking that an equivalent of MA is the topological statement
that, in a compact Hausdorff space whose open sets satisfy the countable
chain condition, the union of less than $2^{\aleph _{0}}$ meagre sets is
meagre [Wei], [Fre1]. This identifies MA as a variant of Baire's Theorem,
and gives it a special role in the investigation of the \textit{additivity
properties} etc. of classical ideals such as the null and meagre sets, for
which see [BartJ].

Given its particular usefulness and origin, MA, termed a \textit{Forcing
Axiom,} inspired the search for further, more powerful, forcing axioms. The
first to occupy centre-stage is the Proper Forcing Axiom, PFA. This is an
extension of $MA(\aleph _{1})$, which draws in more model theory. At the
price of replacing \textit{all} the cardinals $\kappa <\mathfrak{c}$ by
allowing just $\kappa =\aleph _{1}$, PFA relaxes the `ccc' restriction. (In
fact, Todor\v{c}evi\'{c} and Veli\v{c}kovi\'{c} ([Tod], [Vel]) showed that
PFA\ implies that $\mathfrak{c}=\aleph _{2},$ so allowing back in all the,
rather few, cardinals $\kappa <\mathfrak{c.)}$ The relaxation widens access
to the class of \textit{proper} partial orders (below), and so asserts:

\bigskip

PFA:\ \textit{for every partial order }$P$\textit{\ that is proper and any
family} $\mathcal{F}$ \textit{with} $|\mathcal{F}|\leq \aleph _{1}$ \textit{%
of dense subsets of }$P,$\textit{\ there is a filter }$G$\textit{\ in }$P$%
\textit{\ which meets each }$D\in \mathcal{F}$.

\bigskip

The definition of properness refers to the interplay between the whole of
the partial order $P$ and those fragments of $P$ that appear in `suitably
rich' countable structures, as follows. A\ partial order $P$ is \textit{%
proper} if, for any regular uncountable cardinal $\kappa $ and countable
model $\mathcal{M\prec }$ $H\mathcal{(\kappa )}$ (the family of sets
hereditarily of cardinal less than $\kappa $ [Dra, Ch. 3 \S 7]) with $P\in
M: $

\textit{for each }$p\in P\cap M$\textit{\ and each }$q\leq p,$\textit{\
every antichain }$A\in M$\textit{\ contains an element }$r$\textit{\
compatible with }$p.$

(This formulation obviates the need to refer to `maximal antichains'.) The
class of proper partial orders includes both those satisfying ccc (which
preserves cardinality, and cofinality) and those with countable closure
(i.e. guaranteeing a lower bound for any decreasing $\omega $-sequence). A
consistency proof for PFA\ needs use of a supercompact cardinal. See [Bau]
for applications and discussion (especially remarks after his Th. 3.1
concerning the need for a supercompact and its `reflection properties'), and
also [Dev2], and the more recent [Moo]. A wider variant still is SPFA, based
on $\aleph _{1}$-semiproper forcing. The maximal version, known as Martin's
Maximum, MM, was introduced by Foreman, Magidor and Shelah [ForMS], and like
PFA needs a supercompact cardinal for a proof of its consistency. Here the
role of $\omega _{1}$ as $\aleph _{1}$ (in merely prescribing a cardinality
bound) changes in order to create an $\omega _{2}$-chain condition, as we
shall see presently. Prominence is given now to the \textit{stationary }%
subsets of $\omega _{1}$ (defined below), cf. \S 5b; these are the
`non-negligible' subsets in relation to coding, and their definition draws
on some associated `large' sets: the subsets that are closed and unbounded
(cofinal) in $\omega _{1}$, with which we begin. A set $C\subseteq \omega
_{1}$ is \textit{closed} if it contains all its limit points (i.e. $\sup
(C\cap \alpha )\in C\ $for limit $\alpha $ whenever $C\cap \alpha $ is
cofinal in $\alpha );$ such sets form a filter, as any two unbounded closed
sets meet. A subset $S\subseteq \omega _{1}$ is \textit{stationary} if $S$
meets every closed unbounded set. In MM, the partial orders $P$ are required
to preserve stationarity. This condition is motivated by a question about
the \textit{non-stationary ideal, }the ideal of non-stationary sets (denoted 
$\ell _{\text{NS}}$ or NS$_{\omega _{1}}$): whether it is $\omega _{2}$%
-saturated, i.e. whether every $\omega _{2}$-sequence of stationary sets
contains at least two members intersecting again in a stationary set. If so,
then the Boolean algebra $\wp (\omega _{1})/\ell _{\text{NS}}$ is complete
and satisfies the $\omega _{2}$-chain condition. MM implies this.

Woodin [Woo1,2] has forcefully argued for a canonical model where CH\ fails
(cf. Coda); it is a forcing extension of $L(\mathbb{R}),$ i.e. of the Hajnal
`constructible closure' of $\mathbb{R}$ (the class of sets constructible
from some real in $V\ $-- [Dra, Ch. 5 \S 6.1], cf. [Kan, Ch. 1, \S 3]; this
is not to be confused with the L\'{e}vy class of sets `constructible
relative to a given set' [Dra, Ch. 5 \S 6.2], which occurs in \S 5 in the
shape of $L[\mathcal{U}]$ with distinct notation).

\bigskip

\textbf{7. Suslin, Luzin, Sierpi\'{n}ski and their legacy: infinite games
and large cardinals}

After the (necessarily) extensive excursion into logic and model theory, we
now re-anchor all this to analytic practice. Henceforth, we intertwine these
two aspects. For the Analysts's point of view of set theory, we can do no
better at this point than to cite C. A. (Ambrose) Rogers, a modern-day
analyst par excellence (with a pedigree of: Geometry of Numbers, Discrete
geometry, Convexity, Hausdorff measures, Topological descriptive set
theory). In his last phase (post 1960), Rogers famously `would often give
talks entitled \textquotedblleft Which sets do we need?\textquotedblright ,
his answer being: analytic sets' (cited from [Ost7]). To these we now turn.
For background here, see [Rog].

\bigskip

\textbf{7a. Analytic sets}

Analytic subsets of $\mathbb{R}$ are precisely the sets that arise as
projections of planar Borel sets. Their initial (`classical') study,
principally by Suslin, Luzin and Sierpi\'{n}ski, was prompted by Lebesgue's
erroneous assertion, in the course of his research on functions that are
`analytically representable', that these projections were Borel. But they
need not be, as was first observed by Suslin in 1916. Indeed, an analytic
set is Borel iff its complement is also Borel [Sou]. Until that moment the
typical sets considered by analysts were Borel. Fortunately for Lebesgue's
research goals, analytic sets are extremely well-behaved: in the first place
projections of analytic sets are inevitably analytic, and furthermore they
have the following three regularity properties (`the classical regularity
properties' below):\ they are measurable [Lus], they have the property of
Baire [Nik], and likewise the perfect-set property [Ale] (they are either
countable or contain a perfect set), and in certain circumstances are well
approximable from within by compact subsets (they are `capacitable' -- a
property discovered independently by R. O. Davies [Dav] in 1952 and in a
general topological context by G. Choquet in 1952 [Cho1-4]).

The newly discovered sets emerged as the first-level sets of the \textit{%
projective hierarchy} (also called the \textit{Luzin hierarchy}) generated
from the Borel sets by alternately applying the operation of projection and
complementation (a fact later recognized also through the analysis of their
logical complexity:\ counting how many \textit{alternations} of existential
and universal quantifiers over the reals are needed to define them, and
identifying which the \textit{preliminary} quantifier is: existential or
universal). However, the very successful classical study of analytic sets
struggled to promote much of the `good behaviour' up the hierarchy. At the
margins, of particular interest, was K\^{o}ndo's uniformization theorem of
1939 (that a co-analytic planar set has a co-analytic \textit{uniformization}%
, i.e. contains a co-analytic graph selecting one point from each vertical
section).

The message from set theory in G\"{o}del's inner universe of sets $L$ was
particularly depressing: K\^{o}ndo's theorem implied the existence in $L$ of
an analytic sets whose complement failed to have the perfect-set property
(the culprit was the well-ordering of $L$, which relative to $L$ lies at the
second projective level).

Further progress seemed doomed. But an unlikely development, in the shape of
a game-theoretic rival to AC, unblocked the log-jam. However, it was left to
a later generation to pore over the classical achievements to extract the
necessary inspiration from the classicists by drawing in a further theme:
the Banach-Mazur games.

To explain this development we need to explore some analytic-set theory.
Suslin's characterization [Sou] in 1917 of analytic sets $S\subseteq \mathbb{%
R}$ asserts they may be represented in the form%
\[
S=\dbigcup\nolimits_{\mathbf{i}\in \mathbb{N}^{\mathbb{N}}}F(\mathbf{i}%
)=\dbigcup\nolimits_{\mathbf{i}\in \mathbb{N}^{\mathbb{N}}}\dbigcap%
\nolimits_{n=1}^{\infty }F(\mathbf{i}|n), 
\]%
where each of the \textit{determining sets} $F(\mathbf{i}|n)$ is closed and
of diameter at most $2^{-n}$ -- so that $F(\mathbf{i})$ has at most one
member; here%
\[
\mathbf{i}|n:=(i_{0},...,i_{n-1}). 
\]%
(For this reason, the operation taking a determining system to the set $S$
above is now usually called the Suslin operation, though it is sometimes
called the $A$-operation as in [Kur], apparently named for P. S. Alexandrov,
who had devised it to construct perfect subsets of uncountable Borel sets
[Ale].) Implicit in the formula is an operation on the determining system of
sets $\langle F(\mathbf{i}|n):\mathbf{i}|n\in \mathbb{N}^{<\mathbb{N}%
}\rangle ,$ which includes countable intersection and countable union (and
preserves analyticity if the determining system comprises analytic sets
[Rog, Part 1, \S 2.3]). This goes beyond countable union seemingly towards a
continuum union, but one that is constrained by the upper (h)semi-continuity
of the map $\mathbf{i\mapsto }F(\mathbf{i})$.

Under this `continuous union' lie hidden the countable ordinals, by virtue
of the countable tree $T$ of all finite sequences $\mathbf{i}|n$ (ordered by
sequence extension). For any $x$ the associated subtree 
\[
T_{x}:=\{\mathbf{i}|n:x\in F(\mathbf{i}|n)\} 
\]%
is well-founded iff $x\notin S,$ as then $T_{x}$ has no paths (infinite
branches); indeed $x\notin F(\mathbf{i})$ for all $\mathbf{i}$. (This tree
idea, with the $\mathbf{i}|n$ replaced by rationals, goes back, albeit under
the name `sieve', to Lebesgue's construction of a measurable set that is not
Borel.) The overall complexity of the subtree may then be measured by a
countable ordinal, known as the \textit{Luzin-Sierpi\'{n}ski index} of the
tree $T_{x}$ (or of the point $x$) -- [LusS]. This is obtained rather as the
Cantor-Bendixon index of a scattered set is obtained by the repeated
(inductive) removal of isolated points, except that here one removes at each
stage the terminal nodes of a tree. (A moment's reflection shows this
corresponds to a linear ordering of the finite sequences, akin to
lexicographic but adjusted to allow shorter sequences to preceed their
longer extensions, such that the tree is well-ordered iff it is
well-founded: this is the \textit{Kleene-Brouwer order}.)

When the determining system of $S$ (i.e. the family of sets $F(\mathbf{i}|n)$
above) consists of closed sets, it readily follows, via its countable
transfinite definition, that the set of points $x$ in the complement of $S$
with index bounded by $\alpha <\omega _{1}$ is Borel. It is also immediate
that the complement of an analytic set is a union of $\omega _{1}$ Borel
sets, since the index is bounded by $\omega _{1}$. The important \textit{%
boundedness property} of the index (that it remains bounded over any
analytic set $S^{\prime }$ in the complement of $S$ by a corresponding
countable ordinal, a matter that hinges on the `continuous union' aspect)
leads to a proof of the \textit{First} \textit{Separation Theorem}: disjoint
analytic sets may be covered by disjoint Borel sets. From here, as an
immediate corollary, an analytic set with analytic complement is Borel.

\bigskip

\textbf{7b. Banach-Mazur games and the Luzin hierarchy}

We recall that a Banach-Mazur game with target set $S\subseteq \mathbb{R}$
is an infinite positional game which may be viewed as played by two players
`alternately picking ad infinitum' the digits of a decimal expansion of a
real number -- but this needs the interpretation that each player selects a
function (a strategy) determining that player's choice of next digit, given
the current position -- with the first player declared the winner iff the
real number generated from the play of the two strategies falls in $S,$ and
otherwise the second. The target set $S\ $is said to be \textit{determined}
if one or other of the players has a winning strategy. Mazur proposed the
game (this is Problem 43 in the Scottish Book, [Mau]), and Banach responded
in 1935 by characterizing determinacy by the property of Baire. See [GalMS]
for an alternative infinite game which offers a measure-theoretic result as
a contrast to Banach's category result.

It is clear from its description that the game offers a natural
interpretation for a sequence of choices in a manner related to the
countable axiom of choice. In 1962 Mycielski and Steinhaus [MycS] proposed
the Axiom of Determinacy, AD, as an alternative to AC -- in essence setting
the task of ascertaining its consistency relative to ZF. See [Myc] for an
account of the consequences of AD current in 1964, making the case that, in
a hoped-for subuniverse of sets in which AD holds, the well-known
`paradoxes' (Hausdorff, Banach-Tarski, ...) flowing from AC would be ruled
out, while at the same time preserving standard analysis in $\mathbb{R}$
(since `countable choice' for a countable family with union at most a
continuum of members follows from AD -- so, in view of the continuum
restriction, it is usual to work with AD+DC).

We may pass now to a generalization of Suslin's representation for analytic
sets, which enabled higher-level analogues of the classical regularity
properties. Interpreting $\mathbb{N}^{\mathbb{N}}$ as the set of irrationals
(via continued fraction expansion), we may w.l.o.g. assume that $S\subseteq 
\mathbb{N}^{\mathbb{N}}.$ This carries the simplifying advantage that,
ignoring a countable set of lines, we may easily identify planar sets,
regarded as lying in $\mathbb{N}^{\mathbb{N}}\times \mathbb{N}^{\mathbb{N}},$
with subsets of $\mathbb{N}^{\mathbb{N}}$ (merging a pair $(x,y)$ into a
single sequence $\langle x,y\rangle $) and so regard projection as an
operation from $\mathbb{N}^{\mathbb{N}}$ to $\mathbb{N}^{\mathbb{N}}.$

Replacing $F(\mathbf{i}|n)$ by its $2^{-n}$ open swelling $S(\mathbf{i}|n)$
yields that $s\in S$ iff for some $\mathbf{i}\in \mathbb{N}^{\mathbb{N}}$%
\[
s|n\in S(\mathbf{i}|n)\qquad (n\in \mathbb{N)}; 
\]%
here we interpret $s|n$ as a (rational) point of $\mathbb{R}$ (and
implicitly refer to the metric of first difference: $d(x,y)=2^{-n},$ when $%
x,y$ differ first in their $n^{\text{th}}$ term). We can tidy up further
while working in $\mathbb{R},$ by assuming compact $F(\mathbf{i}|n)$ and
replacing $S(\mathbf{i}|n)$ with a union of a finite number of
rational-ended closed intervals. Coding such finite unions by $\mathbb{N}$,
we arrive at a reformulation of Suslin's characterization: for $T\ $a tree
of finite (pairs $(u,v)$ of ) sequences, define the projection of $T$ into $%
\mathbb{N}^{\mathbb{N}}$ by 
\[
p(T):=\{x\in \mathbb{N}^{\mathbb{N}}:(\forall n)[(x|n,\mathbf{i}|n)\in T]\}; 
\]%
then $S\ $is analytic iff $S=p(T)$ for some appropriate tree $T$ of finite
sequences of elements of $\mathbb{N}\times \mathbb{N}$. The generalization
to a $\gamma $\textit{-Suslin} set for ordinals $\gamma $ is obtained by
taking trees $T$ of finite sequences of elements from $\mathbb{N}\times
\gamma ,$ and provides the context allowing the regularity properties of
category and measure to be lifted up the projective hierarchy.

A $\gamma $-Suslin set is said to be a \textit{homogeneously Suslin set }if
there is an $\omega _{1}$-complete ultrafilter $\mathcal{U}_{x|n}$ for each $%
x|n$ such that for all $n$%
\[
\{\mathbf{i}|n\in \gamma ^{n}:(x|n,\mathbf{i}|n)\in T\}\in \mathcal{U}_{x|n} 
\]%
(membership witnessed via a `large' set of nodes), and 
\[
p(T)=\{x\in \mathbb{N}^{\mathbb{N}}:(\forall A_{n}\in \mathcal{U}_{x|n})%
\text{ }\exists \mathbf{i}\text{ }\forall n\text{ }[\mathbf{i}|n\in A_{n}]\} 
\]%
(projection equivalent to passage through a `large' sets of nodes at each
height/level; the sequence $\langle \mathcal{U}_{x|n}:n\in \mathbb{N}\rangle 
$ is then said to be countably complete). In using the index set $\gamma
^{<\omega }$ these generalizations sound muted echoes of the non-separable
theory of analytic sets (pioneered by A. H. Stone and R. W. Hansell -- see
[Sto] and [Ost3]).

Martin, generalizing [Mar], shows in [MarSt, Th. 2.3] that homogeneously $%
\gamma $-Suslin sets are determined (as well as having the classical
regularity properties), and that if Ramsey cardinals exist, then co-analytic
sets are homogenously Souslin. This last result is a re-interpretation of
Martin's earlier theorem [Mar] that if there is a Ramsey cardinal (e.g. if
there is a measurable cardinal), then analytic games are determined.

Two features of the analysis of a co-analytic set $C$ via the Luzin-Sierpi%
\'{n}ski index are of great significance to the study of projective sets.
First, the index maps to the ordinals, i.e. into a well-ordered set, and so
the index induces a prewellordering, rather than a well-ordering on the set $%
C$ (as distinct points of $C$ may be mapped to the same ordinal). Secondly,
denoting the index by $\rho ,$ the relation%
\[
R^{+}(x,y):=x\in C\text{ and }\rho (x)\leq \rho (y), 
\]%
and its negation $R^{-}(x,y)$ are both Borel, and so both co-analytic.
Taking an abstract viewpoint, a class $\Gamma $ of sets in $\mathbb{N}^{%
\mathbb{N}}$ may be said to have the \textit{prewellordering property} if
for every set $C\in \Gamma $ there is a map $\rho :C\rightarrow On$ such
that both of $R^{\pm }(x,y)$ are in $\Gamma .$ (The map is then called a $%
\Gamma $-\textit{norm}.) Suppose that the complementary class $\check{\Gamma}
$ (i.e. of sets with complement in $\Gamma $) is, like the analytic sets,
closed under projection; then the class of sets $\exists ^{1}\Gamma $
obtained as the projections of sets in $\Gamma $ also has the
prewellordering property. This would have been clear to Luzin and Sierpi\'{n}%
ski; but, with the introduction of determinacy, a new feature arises:

\bigskip

\textit{The First Periodicity Theorem} ([Mar], [Mos2]): For a class of sets $%
\Gamma $ for which the sets in the\textit{\ ambiguous class} $\Delta
_{\Gamma }:=\Gamma \cap \check{\Gamma}$ are determined: for every $C\in
\Gamma ,$ if $C\ $admits a $\Gamma $-norm, then $\{y:\forall x[\langle
x,y\rangle \in C]\}$ admits a norm in the class of sets $\forall ^{1}\exists
^{1}\Gamma ,$ i.e. in the class of sets of the form $\forall x\exists
y[\langle x,y\rangle \in C^{\prime }]$ for some $C^{\prime }$ in $\Gamma .$

\bigskip

Thus, in particular: inductively, if the $\mathbf{\Sigma }_{2n}^{1}$-class
(for the $\mathbf{\Sigma }$ and $\mathbf{\Pi }$ notation of the projective
hierarchy, again see \S 9) has the prewellordering property, then so does
the $\mathbf{\Pi }_{2n+1}^{1}$-class, assuming determinacy of the ambiguous
class $\mathbf{\Delta }_{2n}^{1}$. The $\mathbf{\Pi }_{2n+1}^{1}$-class
yields quite directly a prewellordering for the class $\mathbf{\Sigma }%
_{2n+2}^{1}:$ if $A(x)\equiv (\exists y)[\langle x,y\rangle \in C]$ for $C$
in $\mathbf{\Pi }_{2n+1}^{1}$ with norm $\rho _{C},$ then a norm (of the
corresponding class) for $A$ may be defined by 
\[
\rho _{A}(x):=\min \{\rho _{C}(x,y):\langle x,y\rangle \in C\}. 
\]%
Thus the prewellordering property `zig-zags' between the $\mathbf{\Pi }$ and 
$\mathbf{\Sigma }$ classes.

Part of the motivation to take a game-theoretic approach to the projective
sets was the appearance in 1967 of a new proof of the earlier mentioned
Suslin separation theorem for analytic sets given by David Blackwell [Blac]
on the basis of the Gale-Stewart proof of the determinacy of open sets
[GalS] of 1953. The wealth of insights thereafter is history: witness the
very title of Mathias's `Surrealist landscape with figures' survey [Mat],
capturing the spirit of the time.

It was a careful reading of K\^{o}ndo's proof of the uniformization of $%
\mathbf{\Pi }_{1}^{1}$-sets by a $\mathbf{\Pi }_{1}^{1}$ graph that
initially led Moschovakis to isolate a more general kind of $\Gamma $-norm:
that of a $\Gamma $-scale which refers to an $\omega $-sequence of $\Gamma $%
-norms $\rho _{m}$ defined on a set $C$ of $\Gamma $ with associated
relations $R^{\pm }(m,x,y)$ in $\Gamma $ (as with the single $\Gamma $-norm
above), but with an additional `convergence-guiding' property:

for any sequence $c_{n}\in C$ with $c_{n}\rightarrow c_{0}$, if for each $m$ 
\[
\langle \rho _{m}(c_{n}):n\in \omega \rangle \text{ is eventually a
constant, }\lambda _{m}\text{ say,} 
\]%
then $c_{0}\in C$ and $\rho _{m}(c_{0})\leq \lambda _{m}$ for all $m.$ (See
e.g. [MarK, \S 8.2].)

Mutatis mutandis, the Moschovakis \textit{Second Periodicity Theorem} [Mos2]
has the same form as the First but with $\Gamma $-scale replacing $\Gamma $%
-norm throughout. Analogously, the Second Theorem implies that the K\^{o}ndo
uniformization property likewise zigzags between the $\mathbf{\Pi }$ and $%
\mathbf{\Sigma }$ classes -- see [Mos1].

Guided by the original $\mathbf{\Pi }_{1}^{1}$-norm (the Luzin-Sierpi\'{n}%
ski index), having range in $\omega _{1}$ (less, if the $\mathbf{\Pi }%
_{1}^{1}$ set in question is Borel), one defines the \textit{projective
ordinal} of level $n$ by reference to the sets in the ambiguous class $%
\mathbf{\Delta }_{n}^{1}$%
\[
\mathbf{\delta }_{n}^{1}:=\text{supremum of the lengths of prewellorderings
in }\mathbf{\Delta }_{n}^{1}. 
\]%
Martin showed that $\mathbf{\delta }_{2}^{1}\leq \omega _{2},$ with equality
implied under AD by the Moschovakis result that $\mathbf{\delta }_{n}^{1}$
for $n\geq 1$ is a cardinal and that, under PD, $\mathbf{\delta }_{2n}^{1}<%
\mathbf{\delta }_{2n+2}^{1}.$ Under AD+DC $\mathbf{\delta }_{2n}^{1}=(%
\mathbf{\delta }_{2n-1}^{1})^{+}$ (i.e. the even-indexed ordinal is the
successor of the preceding odd-indexed one); furthermore, Jackson's theorem
[Jac1,2] asserts that under AD+DC%
\[
\mathbf{\delta }_{2n-1}^{1}=\aleph _{w(2n-1)+1}, 
\]%
where $w$ is defined via iterated exponentiation: inductively by $%
w(m+1)=\omega ^{w(m)}$ with $w(1)=\omega .$

A concerted effort to assess the consistency strength of the determinacy
assumption for $\mathbf{\Delta }_{n}^{1}$ ultimately led to the result that
this is equiconsistent with the existence of $n$ Woodin cardinals below a
measurable cardinal.

\bigskip

\textbf{8. Shadows}

Here we wrap up our survey of the set-theoretical domain. We have seen how
combinatorial properties, some `high up', in Cantor's world affect
properties of the real line down below. When powerful axioms extend familiar
properties in desirable ways one is led to ask whether one can get away with
less and get if not the same outcome, then `almost' the same (in some
sense). To this end Mycielski and Tomkowicz [MycT] speak in very suggestive
language of \textit{shadows} of AC in their chosen setting of $L(\mathbb{R}%
), $ a model of set theory that resolves some of the hardest set-theory
problems. Their quest is theorems of ZFC\ that have corollaries that are
theorems of ZF+AD -- see [MycT]. In $L(\mathbb{R})$ AD\ implies DC [Kec1],
and the present authors have come to view DC as a natural ally for analysis.
We give our favourite example of this, and then, after a brief review of
syntactical teminology in \S 9, we survey in \S 10 results which give
further succour, if one is willing in the interests of plurality to conduct
mathematics in an appropriate helpful (indeed playful, to borrow the term
from [Mos1], when games are enlisted) subuniverse.

\textit{An example with the Principle of D\textit{ependent Choice DC in mind.%
} }We begin with an example concerned with real-valued sublinear functions
on $\mathbb{R}$ which `almost' follow Banach's enduring paradigmatic
definition. They are subadditive, i.e. satisfying $f(x+y)\leq f(x)+f(y),$
but in one variant they are $\mathbb{N}$-homogeneous in the sense that $%
f(nx)=nf(x)$ for $n=0,1,2,...$, so $\mathbb{Q}_{+}$-homogeneous, and for all 
$x.$ In other variants the quantification over $x$ may also be thinned --
see [BinO6]. In electing to study sublinear functions as possible
realizations of norms, Berz ([Berz], [BinO6]) showed, for measurable $f,$
that the graph of $f$ is conical -- comprises two half lines through the
origin; however, his argument relied on AC, in the usual form of \textit{%
Zorn's Lemma}, which he used in the context of $\mathbb{R}$ over the field
of scalars $\mathbb{Q}$. In spirit he follows Hamel's construction of a
discontinuous additive function, and so ultimately this rests on transfinite
induction of \textit{continuum} length requiring \textit{continuum} many
selections. Our own proof [BinO6] (cf. [BinO7, 10]) of Berz's theorem, taken
in a wider context including Banach spaces, depends in effect on the Baire
Category Theorem BC, or the completeness of $\mathbb{R}$ (in either of the
distinct roles of `Cauchy-sequential' and `Cauchy-filter' completeness, the
latter stronger in the absence of AC, see [FosM, \S 3] and also [DodM, \S 7, 
\S 2]): we rely on generalizations of the \textit{Kestelman-Borwein-Ditor
Theorem}, KBD, asserting that for any (category/measure theoretic)
non-negligible set $T\ $and any null sequence $z_{n}\rightarrow 0,$ for
quasi all $t\in T$ the $t$-translate of some subsequence $z_{n(m)}$
(dependent on $t$) embeds in $T,$ i.e. $t+z_{n(m)}\in T.$ See [MilO] for a
discussion of this `shift-compactness' notion. KBD is a variant of BC. So
the proof ultimately rests on elementary induction via the \textit{Axiom
(Principle) of Dependent Choice(s)} DC (thus named in 1948 by Tarski [Tar2,
p. 96] and studied in [Most], but anticipated in 1942 by Bernays [Ber, Axiom
IV*, p. 86] -- see [Jec1, \S 8.1], [Jec2, Ch. 5]); DC in turn is equivalent
to BC by a result of Blair [Bla]. (For further results in this direction see
also [Pin3,4], [Gol], [HerK], [Wol], and the textbook [Her].)

The relevance of KBD in the setting of a Polish group comes from its various
corollaries which include the Steinhaus-Weil Interior points Theorem
[BinO9], the Open Mapping Theorem and its generalization to group actions:
the Effros Theorem -- see [vMil], [Ost4,5,6]. For a target set $T$ that is a
dense $G_{\delta }$, embeddings which are performed simultaneously in any
neighbourhood by a perfect subset of $T$ of a fixed set $Z$ (not necessarily
a null sequence) into $T$ characterize those sets $Z$ that are strong
measure zero -- see[GalMS].

We note that DC is equivalent to a statement about trees: a pruned tree has
an infinite branch (for which see [Kec2, 20.B]) and so by its very nature is
an ingredient in set-theory axiom systems which consider the extent to which
Banach-Mazur-type games (with underlying tree structure) are determined. The
latter in turn have been viewed as generalizations of Baire's Theorem ever
since Choquet [Cho5] -- cf. [Kec2, 8C,D,E]. Inevitably, determinacy and the
study of the relationship between category and measure go hand in hand.

\bigskip

\textbf{9. The syntax of Analysis: Category/measure regularity versus
practicality}

The Baire/measurable property discussed at various points above is usually
satisfied in mathematical practice. Indeed, any analytic subset of $\mathbb{R%
}$ possesses these properties ([Rog, Part 1 \S 2.9], [Kec2, 29.5]), hence so
do all the sets in the $\sigma $-algebra that they generate (the $C$-sets,
[Kec2, \S 29.D], $C$ for \textit{cribl\'{e} }-- see [Bur1, 2], cf. [BinO4]).
There is a broader class still. Recall first that an analytic set may be
viewed as a projection of a planar Borel set $P,$ so is definable as $%
\{x:\Phi (x)\}$ via the $\mathbf{\Sigma }_{1}^{1}$ formula $\Phi
(x):=(\exists y\in \mathbb{R})[(x,y)\in P];$ here the notation $\mathbf{%
\Sigma }_{1}^{1}$ indicates one quantifier block (the subscripted value) of
existential quantification, ranging over reals (type 1 objects -- the
superscripted value). Use of the \textit{bold-face version} of the symbol
indicates the need to refer to \textit{arbitrary} coding (by reals not
necessarily in an \textit{effective} manner, for which see [Gao, \S 1.5]) of
the various open sets needed to construct $P.$ (An open set $U$ is \textit{%
coded} by the sequence of rational intervals contained in $U.)$ Effective
variants are rendered in light-face.

Consider a set $A$ such that both $A$ and $\mathbb{R}\backslash A$ may be
defined by a $\mathbf{\Sigma }_{2}^{1}$ formula, say respectively as $%
\{x:\Phi (x)\}$ and $\{x:\Psi (x)\}$, where $\Phi (x):=(\exists y\in \mathbb{%
R})(\forall z\in \mathbb{R})(x,y,z)\in P\}$ now, and similarly $\Psi .$ This
means that $A$ is both $\mathbf{\Sigma }_{2}^{1}$ and $\mathbf{\Pi }_{2}^{1}$
(with $\mathbf{\Pi }$ indicating a leading universal quantifier block), and
so is in the ambiguous class $\mathbf{\Delta }_{2}^{1}.$ If in addition the
equivalence%
\[
\Phi (x)\Longleftrightarrow \lnot \Psi (x) 
\]%
is provable in $ZF,$ i.e. \textit{without reference} to $AC,$ then $A$ is
said to be \textit{provably} $\mathbf{\Delta }_{2}^{1}.$ It turns out that
such sets have the Baire/measurable property -- see [FenN], where these are
generalized to the \textit{universally (=absolutely) measurable }sets (cf.
[BinO6, \S 2]); the idea is ascribed to Solovay in [Kan, Ch. 3 Ex. 14.4].
How much further this may go depends on what axioms of set theory are
admitted, a matter to which we now turn.

Our interest in such matters derives from the \textit{Character Theorems} of
regular variation, noted in [BinO3, \S 3] (revisited in [BinO5, \S 11]),
which identify the logical complexity of the function 
\[
h^{\ast }(x):=\lim \sup_{t\rightarrow \infty }h(t+x)-h(t), 
\]%
which is $\mathbf{\Delta }_{2}^{1}$ if the function $h$ is Borel (and is $%
\mathbf{\Pi }_{2}^{1}$ if $h$ is analytic, and $\mathbf{\Pi }_{3}^{1}$ if $h$
is co-analytic). We argued in [BinO3, \S 5] that $\mathbf{\Delta }_{2}^{1}$
is a natural setting in which to study regular variation.

\bigskip

\textbf{10. Category-Measure duality}

\noindent \textit{1. Practical axiomatic alternatives: LM, PB, AD, PD. }%
While ZF is common ground in mathematics, AC is not, and alternatives to it
are widely used, in which for example all sets are Lebesgue-measurable
(usually abbreviated to LM) and all sets have the Baire property, sometimes
abbreviated to PB (as distinct from BP to indicate individual `possession of
the Baire property'). One such is DC above. As Solovay [Sol3, p. 25] points
out, this axiom is sufficient for the establishment of Lebesgue measure,
i.e. including its translation invariance and countable additivity
("...positive results ... of measure theory..."), and may be assumed
together with LM. Another is the \textit{Axiom of Determinacy} AD mentioned
above and introduced by Mycielski and Steinhaus [MycS]; this implies LM, for
which see [MycSw], and PB, the latter a result, mentioned in \S 7, due to
Banach -- see [Kec2, 38.B]. Its introduction inspired remarkable and still
current developments in set theory concerned with determinacy of `definable'
sets of reals (see [ForK] and particularly [Nee]) and consequent
combinatorial properties (such as partition relations) of the alephs (see
[Kle]); again see \S 7. Others include the (weaker) \textit{Axiom of
Projective Determinacy} PD [Kec2, \S\ 38.B], cf. \S 7, restricting the
operation of AD\ to the smaller class of projective sets. (The independence
and consistency of DC\ versus AD was established respectively in Solovay
[Sol4] and Kechris [Kec2] -- see also [KecS]; cf. [DalW1], [Ost2].)

\noindent \textit{2. LM versus PB. }In 1983 Raissonier and Stern [RaiS, Th.
2] (cf. [Bart1,2]), inspired by then current work of Shelah (circulating in
manuscript since 1980) and earlier work of Solovay, showed that \textit{if
every }$\Sigma _{2}^{1}$\textit{\ set is Lebesgue measurable, then every }$%
\Sigma _{2}^{1}$\textit{\ set has BP}, whereas the converse fails -- for the
latter see [Ster] -- cf. [BartJ, \S 9.3] and [Paw]. This demonstrates that
measurability is in fact the stronger notion -- see [JudSh, \S 1] for a
discussion of the consistency of analogues at level 3 and beyond -- which is
one reason why we regard category rather than measure as \textit{primary}.
For we have seen above how the category version of Berz's theorem implies
its measure version; see also [BinO6, 10].

Note that the assumption of G\"{o}del's \textit{Axiom of Constructibility} $%
V=L,$ a strengthening of AC, yields $\mathbf{\Delta }_{2}^{1}$
non-measurable subsets, so that the Fenstad-Normann result on the narrower
class of provably $\mathbf{\Delta }_{2}^{1}$ sets mentioned in \S 9 above
marks the limit of such results in a purely ZF framework (at level 2).

\noindent \textit{3. Consistency and the role of large cardinals. }While LM\
and PB are inconsistent with AC, such axioms can be consistent with DC.
Justification with scant exception involves some form of large-cardinal
assumption, which in turn, as in \S 4, calibrates relative consistency
strengths -- see [Kan] and [KoeW] (cf. [Lar] and [KanM]). Thus Solovay
[Sol3] in 1970 was the first to show the equiconsistency of ZF+DC+LM+PB with
that of ZFC+`\textit{there exists an inaccessible cardinal}'. The appearance
of the inaccessible in this result is not altogether incongruous, given its
emergence in results (from 1930 onwards) due to Banach [Ban] (under GCH),
Ulam [Ula] (under AC), and Tarski [Tar1], concerning the cardinalities of
sets supporting a countably additive/finitely additive [0,1]-valued/$\{0,1\}$%
-valued measure (cf. [Bog, 1.12(x)], [Fre2]). Later in 1984 Shelah [She1,
5.1] showed in ZF+DC that already the measurability of all $\mathbf{\Sigma }%
_{1}^{3}$ sets implies that $\aleph _{1}^{L}$ is inaccessible (the symbol $%
\aleph _{1}^{L}$ refers to the substructure/ subuniverse of constructible
sets and denotes the first uncountable ordinal therein -- cf. \S 2). As a
consequence, Shelah [She1,5.1A] showed that ZF+DC+LM\ is equiconsistent with
ZF+`\textit{there exists an inaccessible}', whereas [She1, 7.17] ZF+DC+PB\
is equiconsistent with just ZFC (i.e. without reference to inaccessible
cardinals), so driving another wedge between classical measure-category
symmetries (see [JudSh] for further, related `wedges'). The latter
consistency theorem relies on the result [She1, 7.16] that any model of ZFC\
+ CH has a generic (forcing) extension satisfying ZF+ `\textit{every set of
reals (first-order) defined using a real and an ordinal parameter has BP}'.
(Here `first-order' restricts the range of any quantifiers.) For a
topological proof see Stern [Ster].

\noindent \textit{4. LM versus PB continued. }Raisonnier [Rai, Th. 5] (cf.
[She1, 5.1B]) has shown that in ZF+DC one can prove that if there is an
uncountable well-ordered set of reals (in particular a subset of cardinality 
$\aleph _{1}$), then there is a non-measurable set of reals. (This motivates
Judah and Spinas [JudSp] to consider generalizations including the
consistency of the $\omega _{1}$-variant of DC.) See also Judah and Ros\l %
anowski [JudR] for a model (due to Shelah) in which ZF+DC+LM+$\lnot $PB\
holds, and also [She2] where an inaccessible cardinal is used to show
consistency of ZF+LM+$\lnot $PB+`\textit{there is an uncountable set without
a perfect subset}'. For a textbook treatment of much of this material see
again [BartJ].

Raisonnier [Rai, Th. 3] notes the result, due to Shelah and Stern, that
there is a model for ZF+DC+PB+$\aleph _{1}=\aleph _{1}^{L}$+ `\textit{the
ordinally definable subsets of reals are measurable}'. So, in particular by
Raisonnier's result, there is a non-measurable set in this model. Shelah's
result indicates that the non-measurable is either $\Sigma _{3}^{1}$
(light-face symbol: all open sets coded effectively) or $\mathbf{\Sigma }%
_{2}^{1}$ (bold-face). Thus here PB+$\lnot $LM holds.

\noindent \textit{5. Regularity of reasonably definable sets. }From the
existence of suitably large cardinals flows a most remarkable result due to
Shelah and Woodin [SheW] justifying the opening practical remark about BP,
which is that every `reasonably definable' set of reals is Lebesgue
measurable: compare the commentary in [BecK] following their Th 5.3.2. This
is a latter-day sweeping generalization of a theorem due to Solovay (cf.
[Sol2]) that, subject to large-cardinal assumptions, $\mathbf{\Sigma }%
_{2}^{1}$ \textit{sets are measurable} (and so also \textit{have BP} by
[RaiS]).

\bigskip

\noindent \textbf{Coda}

To return to the algebraic characterization of the reals as `the' complete
archimedean ordered field: it is the `complete' which hides the `modulo
cardinality' and `modulo which sets are available' aspects. It is always
good to look at familiar mathematics, and ask oneself the analogous question
in that context.

As working analysts ourselves, we feel for those of our colleagues new to
these matters, who may look fondly back to an age of `bygone innocence',
when `one didn't need to worry about such things'. We prefer instead to
marvel at the unfathomable richness of mathematics. As usual, Shakespeare
puts his finger on it somewhere:

\begin{center}
{\small \textit{There are more things in heaven and earth, Horatio,}}\\[0pt]
{\small \textit{Than are dreamt of in our philosophy.}}
\end{center}

\noindent So we have only mathematical `gut-feeling and belief', as with
Mickiewicz:

\begin{center}
{\small \textit{Czucie i wiara silniej m\'{o}wi do mnie},}\\[0pt]
{\small \textit{Ni\.{z} m\k{e}drca szkie\l ko i oko.}}\\[0pt]
\end{center}

\noindent -- `Feeling and faith more forcefully persuade, Than the lens and
the eye of a sage'.

Thus it is that we close with two `high-profile' attitudes towards Solovay's
dictum that the continuum `can be anything it ought to be', to both of which
Woodin has contributed. On the one hand there is a putative $L$-like
`ultimate inner model' (leading to $V=Ult$-$L$) [Woo3], which permits
adjunction of known large-cardinal axioms; under it the continuum is $\aleph
_{1}.$ On the other hand is the argument, offered by Woodin in [Woo2], close
in spirit to the Forcing Axioms of \S 8 as it depends on closure under (set)
forcing in the presence of large cardinals; under this the continuum is $%
\aleph _{2}$.

\newpage

\noindent \textbf{References}

\noindent \lbrack Ale] P. S. Aleksandrov, Sur la puissance des ensembles
mesurable B. \textsl{C. R. Acad. Sci. Paris} \textbf{162}, (1916), 323-325.%
\newline
\noindent \lbrack Ban] S. Banach, \textsl{Th\'{e}orie des op\'{e}rations lin%
\'{e}aires.} Monografje Mat. \textbf{1}, Warszawa, 1932; Chelsea Publishing
Co., New York, 1955; in \textsl{Oeuvres}, \textbf{II}, PWN-1979 (see
http://matwbn-old.icm.edu.pl/), \'{E}ditions Jacques Gabay, Sceaux, 1993; 
\textsl{Theory of linear operations.} Translated from the French by F.
Jellett. North-Holland Mathematical Library\textbf{\ 38}, 1987.\newline
\noindent \lbrack Bart1] T. Bartoszy\'{n}ski, Additivity of measure implies
additivity of category. \textsl{Trans. Amer. Math. Soc.} \textbf{281}
(1984), no. 1, 209--213 \newline
\noindent \lbrack Bart2] T. Bartoszy\'{n}ski, \textsl{Invariants of measure
and category}, Ch. 7 in [ForK].\newline
\noindent \lbrack BartJ] T. Bartoszy\'{n}ski and H. Judah, \textsl{Set
Theory: On the structure of the real line}, Peters 1995.\newline
\noindent \lbrack Bau] J. E. Baumgartner, Applications of the Proper Forcing
Axiom, \textsl{Handbook of set-theoretic topology}, (eds K. Kunen and J.E.
Vaughan), 913-960, North-Holland, 1984.\newline
\noindent \lbrack BecK] H. Becker, A. S. Kechris, \textsl{The descriptive
set theory of Polish group actions.} London Math. Soc. Lect. Note Ser., 
\textbf{232}, Cambridge University Press, 1996.\newline
\noindent \lbrack BelS] J. Bell, A. Slomson, \textsl{Models and
ultraproducts: an introduction}, N. Holland, 1969.\newline
\noindent \lbrack Ber] P. Bernays, A system of axiomatic set theory. III.
Infinity and enumerability. Analysis, \textsl{J. Symbolic Logic} 7 (1942),
65--89.\newline
\noindent \lbrack Berz] E. Berz, Sublinear functions on $\mathbb{R}.$ 
\textrm{\textsl{Aequat. Math.}} \textbf{12} (1975), 200-206.\newline
\noindent \lbrack Bin] N. H. Bingham, Finite additivity versus countable
additivity. \textsl{J. Electron. Hist. Probab. Stat.} \textbf{6} (2010), no.
1, 35 pp.\newline
\noindent \lbrack BinGT] N. H. Bingham, C. M. Goldie and J. L. Teugels, 
\textsl{Regular variation}, 2nd ed., Cambridge University Press, 1989 (1st
ed. 1987).\newline
\noindent \lbrack BinO1] N. H. Bingham and A. J. Ostaszewski, Infinite
combinatorics and the foundations of regular variation, \textsl{Journal of
Math. Anal. Appl.} \textbf{360} (2009), 518-529.\newline
\noindent \lbrack BinO2] N. H. Bingham and A. J. Ostaszewski, Beyond
Lebesgue and Baire II: Bitopology and measure-category duality. \textsl{%
Colloquium Math.} \textbf{121} (2010), 225-238.\newline
\noindent \lbrack BinO3] N. H. Bingham and A. J. Ostaszewski, Regular
variation without limits. \textsl{J. Math. Anal. Appl.}, \textbf{370}
(2010), 322-338.\newline
\noindent \lbrack BinO4] N. H. Bingham and A. J. Ostaszewski, Dichotomy and
infinite combinatorics: the theorems of Steinhaus and Ostrowski. \textsl{%
Math. Proc. Camb. Phil. Soc.} \textbf{150} (2011), 1-22. \newline
\noindent \lbrack BinO5] N. H. Bingham and A. J. Ostaszewski, Beurling
moving averages and approximate homomorphisms, \textsl{Indag. Math. (NS) } 
\textbf{27} (2016), 601--633 (fuller version: arXiv1407.4093v2).\newline
\noindent \lbrack BinO6] N. H. Bingham and A. J. Ostaszewski,
Category-measure duality: convexity, mid-point convexity and Berz
sublinearity, \textsl{Aequationes Math.}, \textbf{91.5} (2017), 801--836 (
fuller version: arXiv1607.05750).\newline
\noindent \lbrack BinO7] N. H. Bingham and A. J. Ostaszewski, Additivity,
subadditivity and linearity: Automatic continuity and quantifier weakening,
Indagationes Math., Online, 2017 (arXiv1405.3948v3).\newline
\noindent \lbrack BinO8] N. H. Bingham and A. J. Ostaszewski, Beyond
Lebesgue and Baire IV: Density topologies and a converse Steinhaus-Weil
Theorem, \textsl{Topol. Appl.} (2018)
https://doi.org/10.1016/j.topol.2017.12.029, to appear (arXiv1607.00031).%
\newline
\noindent \lbrack BinO9] N. H. Bingham and A. J. Ostaszewski, The
Steinhaus-Weil property: its converse, Solecki amenability and
subcontinuity, arXiv1607.00049.\newline
\noindent \lbrack BinO10] N. H. Bingham and A. J. Ostaszewski, Variants on
the Berz sublinearity theorem, arXiv1712.05183.\newline
\noindent \lbrack Blac] D. Blackwell, Infinite games and analytic sets, 
\textsl{Proc. Nat. Acad. Sci.} \textbf{58} (1967), 1836-1837

\noindent \lbrack Bla] C. E. Blair, The Baire category theorem implies the
principle of dependent choices. \textsl{Bull. Acad. Polon. Sci. S\'{e}r.
Sci. Math. Astronom. Phys.} \textbf{25} (1977), 933--934.\newline
\noindent \lbrack Bog] V. I. Bogachev, \textsl{Measure theory.} Vol. I.
Springer, 2007. \newline
\noindent \lbrack Bro1] L. E. J. Brouwer, \"{U}ber Abbildung von
Mannigfaltigkeiten. \textsl{Math. Ann.} \textbf{71} (1911), 97--115. \newline
\noindent \lbrack Bro2] L. E. J. Brouwer, An intuitionist correction of the
fixed-point theorem on the sphere. \textsl{Proc. Roy. Soc. London. Ser. A.}%
\textbf{\ 213}, (1952), 1--2.\newline
\noindent \lbrack Bur1] J. P. Burgess, Classical hierarchies from a modern
standpoint. I. C-sets. \textsl{Fund. Math.} \textbf{115 }(1983), 81--95.%
\newline
\noindent \lbrack Bur2] J. P. Burgess, Classical hierarchies from a modern
standpoint. II. R-sets. \textsl{Fund. Math.} \textbf{115 }(1983), 97--105.%
\newline
\noindent \lbrack Car1] H. Cartan, Th\'{e}orie des filtres,\textsl{\ Comptes
Rendus de l'Acad. des Sci.}, Paris, \textbf{205}, (1937), 595-598.\newline
\noindent \lbrack Car2] H. Cartan, Filtres et ultrafiltres, ibid., 205,
777-779, 1937.\newline
\noindent \lbrack Cho1] G. Choquet, Capacit\'{e}s. Premi\`{e}res d\'{e}%
finitions. (French) \textsl{C. R. Acad. Sci. Paris} \textbf{234}, (1952).
35--37.\newline
\noindent \lbrack Cho2] G. Choquet,Extension et restriction d'une capacit%
\'{e}. (French) \textsl{C. R. Acad. Sci. Paris} \textbf{234}, (1952).
383--385.\newline
\noindent \lbrack Cho3] G. Choquet,Propri\'{e}t\'{e}s fonctionnelles des
capacit\'{e}s altern\'{e}es ou monotones. Exemples. (French) C. R. \textsl{%
C. R. Acad. Sci. Paris} \textbf{234}, (1952). 498--500. \newline
\noindent \lbrack Cho4] G. Choquet, Theory of capacities. \textsl{Ann. Inst.
Fourier, Grenoble} \textbf{5} (1953--1954), 131--295 (1955).\newline
\noindent \lbrack Cho5] G. Choquet, \textsl{Lectures on analysis}, Vol. I,
Benjamin, New York, 1969.\newline
\noindent \lbrack Cie] K. Ciesielski, \textsl{Set theory for the working
mathematician.} London Mathematical Society Student Texts \textbf{39}.
Cambridge University Press, 1997.\newline
\noindent \lbrack Coh1] P. Cohen, The independence of the continuum
hypothesis. \textsl{Proc. Nat. Acad. Sci. }\textbf{50} (1963), 1143--1148.%
\newline
\noindent \lbrack Coh2] P. Cohen, The independence of the continuum
hypothesis II. \textsl{Proc. Nat. Acad. Sci. }\textbf{51} (1964), 105-110.%
\newline
\noindent \lbrack Coh3] P. J. Cohen, \textsl{Set theory and the continuum
hypothesis. }W. A. Benjamin, 1966.\newline
\noindent \lbrack Cohn] P. M. Cohn, P. \textsl{Algebra. } Vol. 2. 2$^{\text{%
nd}}$ ed. Wiley 1989.(1$^{\text{st}}$ ed. 1977)\newline
\noindent \lbrack ComN] W. W. Comfort, S. Negrepontis, \textsl{Theory of
ultrafilters}, Grund. math. Wiss. \textbf{211}, Springer, 1974.\newline
\noindent \lbrack DalW1] H. G. Dales and W. H.Woodin, \textsl{An
introduction to independence for analysts.} London Mathematical Society
Lecture Note Series, 115. Cambridge University Press, Cambridge, 1987.
xiv+241 pp.\newline
\noindent \lbrack DalW2] H. G. Dales and W. H.Woodin, \textsl{Super-real
fields. Totally ordered fields with additional structure.} London
Mathematical Society Monographs, New Series, 17, Oxford University Press,
1996.\newline
\noindent \lbrack Dav] R. O. Davies, Subsets of finite measure in analytic
sets. \textsl{Indagationes Math.} \textbf{14}, (1952), 488--489.\newline
\noindent \lbrack Dev1] K. Devlin,\textsl{\ Aspects of constructibility}.
Lecture Notes in Mathematics \textbf{354}, Springer, 1973.

\noindent \lbrack Dev2] K. Devlin, The Yorkshireman's guide to proper
forcing. \textsl{Surveys in set theory}, 60--115, London Math. Soc. Lecture
Note Ser. \textbf{87}, Cambridge Univ. Press, 1983.\newline
\noindent \lbrack DodM] J. Dodu, M. Morillon, The Hahn-Banach property and
the axiom of choice. \textsl{Math. Log. Q. }\textbf{45.3} (1999), 299--314.%
\newline
\noindent \lbrack Dra] F. Drake, \textsl{Set theory: an introduction to
large cardinals}, North-Holland, 1974.\newline
\noindent \lbrack EhrM] A. Ehrenfeucht, A. Mostowski, Models of axiomatic
theories admitting automorphisms. \textsl{Fund. Math.} \textbf{43} (1956),
50--68.\newline
\noindent \lbrack ErdR] P. Erd\H{o}s, R. Rado, A partition calculus in set
theory. \textsl{Bull. Amer. Math. Soc.} \textbf{62} (1956), 427--489.\newline
\noindent \lbrack FalGOS] K. Falconer, P. M. Gruber, A. Ostaszewski and T.
Stuart, Claude Ambrose Rogers, Biogr. Mems Fell. R. Soc. \textbf{61} (2015),
403-435.\newline
\noindent \lbrack FenN] J. E. Fenstad, D. Normann, On absolutely measurable
sets. \textsl{Fund. Math.} \textbf{81.2} (1973/74), 91--98.\newline
\noindent \lbrack ForK] M. Foreman, A. Kanamori, \textsl{Handbook of Set
theory}, Springer, 2010.\newline
\noindent \lbrack ForMS] M. Foreman, M. Magidor, S. Shelah, Martin's
maximum, saturated ideals and nonregular ultrafilters. I. \textsl{Ann. of
Math.} (2) \textbf{127} (1988), no. 1, 1--47.\newline
\noindent \lbrack FosM] J. Fossy, M. Morillon, The Baire category property
and some notions of compactness, \textsl{J. London Math. Soc.} (2) 57
(1998), 1-19.\newline
\noindent \lbrack Fre1] D. Fremlin, \textsl{Consequences of Martin's axiom.}
Cambridge Tracts in Math. \textbf{84}, Cambridge University Press, 1984.%
\newline
\noindent \lbrack Fre2] D. Fremlin, \textsl{Measure theory Vol. 5:
Set-theoretic measure theory, Parts I, II.} Torres-Fremlin, 2008.\newline
\noindent \lbrack GabW] D. M. Gabbay, J. Woods (ed.) \textsl{Handbook of the
History of Logic} vol. \textbf{5}, Logic from Russell to Church,
North-Holland, 2009.\newline
\noindent \lbrack GalS] D. Gale F. M. Stewart, Infinite games with perfect
information, \textsl{Ann. Math. Stud.} \textbf{28} (1953), 245-266.\newline
\noindent \lbrack GalMS] F. Galvin, J. Mycielski, R. M. Solovay, Strong
measure zero and infinite games. \textsl{Arch. Math. Logic} \textbf{56}
(2017), 725--732.\newline
\noindent \lbrack Gao] S. Gao, \textsl{Invariant descriptive theory}, CRC\
Press, 2009.\newline
\noindent \lbrack GarP] R. J. Gardner, W. F. Pfeffer, Borel measures, 
\textsl{Handbook of set-theoretic topology}, (eds K. Kunen and J.E.
Vaughan), 887--911, North-Holland, 1984.\newline
\noindent \lbrack Gol] R. Goldblatt, On the role of the Baire category
theorem and dependent choice in the foundations of logic. \textsl{J.
Symbolic Logic} \textbf{50} (1985), no. 2, 412--422.\newline
\noindent \lbrack Hall1] M. Hallett, Absoluteness and the Skolem paradox.
Logic, mathematics, philosophy: vintage enthusiasms,\textsl{\ West. Ont.
Ser. Philos. Sci.},\textbf{\ 75}, (2011), 189--218.\newline
\noindent \lbrack Hall2] M. Hallett, \textsl{Cantorian set theory and
limitation of size}. Oxford Logic Guides \textbf{10}, Oxford University
Press, 1984.\newline
\noindent \lbrack Hal] P. R. Halmos, \textsl{Naive set theory}. Undergrad.
Texts in Math. Springer, 1974 (reprint of the 1960 edition).\newline
\noindent \lbrack Har] L. Harrington, Analytic determinacy and 0%
. \textsl{J. Symbolic Logic} \textbf{43.4} (1978), 685--693.\newline
\noindent \lbrack Her] H. Herrlich, \textsl{Axiom of choice.} Lecture Notes
in Mathematics, 1876. Springer, 2006.\newline
\noindent \lbrack HerK] H. Herrlich, K. Keremedis, The Baire category
theorem, and the axiom of dependent choice.\textsl{\ Comment. Math. Univ.
Carolin.} \textbf{40} (1999), no. 4, 771-775.\newline
\noindent \lbrack Hew] E. Hewitt, Rings of real-valued continuous functions.
I. \textsl{Trans. Amer. Math. Soc.} \textbf{64}, (1948), 45--99.\newline
\noindent \lbrack Hil1] D. Hilbert, \textsl{Grundlagen der Geometrie.}
Teubner, 1899.\newline
\noindent \lbrack Hil2] D. Hilbert, Les principes fondamentaux de la g\'{e}om%
\'{e}trie. \textsl{Ann. Sci. \'{E}cole Norm. Sup.} (3) \textbf{17} (1900),
103--209.\newline
\noindent \lbrack Hil3] D. Hilbert, Ueber die Grundlagen der Geometrie. 
\textsl{Math. Ann.} \textbf{56} (1902), no. 3, 381--422. \newline
\noindent \lbrack HilB] D. Hilbert, P. Bernays, Grundlagen der Mathematik,
Springer, 1934.\newline
\noindent \lbrack Hod] W. Hodges, A shorter model theory, Cambridge
University Press, 1997.\newline
\noindent \lbrack HowR] P. Howard, J. E. Rubin, The Boolean prime ideal
theorem plus countable choice do [does] not imply dependent choice. \textsl{%
Math. Logic Quart.} \textbf{42} (1996), no. 3, 410--420.\newline
\noindent \lbrack Jac1] S. Jackson, AD and the projective ordinals. \textsl{%
Cabal Seminar} 81--85, 117--220, Lecture Notes in Math. \textbf{1333},
Springer, 1988.\newline
\noindent \lbrack Jac2] S. Jackson, \textsl{Structural consequences of AD},
Ch. 21 in [ForK].\newline
\noindent \lbrack Jec1] T. J. Jech, \textsl{The axiom of choice.} Studies in
Logic and the Foundations of Mathematics, Vol. \textbf{75}. North-Holland,
1973.\newline
\noindent \lbrack Jec2] T. J. Jech, \textsl{Set Theory}, 3$^{\text{rd}}$
Millennium ed. Springer, 2003.\newline
\noindent \lbrack JudR] H. Judah, A. Ros\l anowski, On Shelah's
amalgamation. \textsl{Set theory of the reals} (Ramat Gan, 1991), 385--414,
Israel Math. Conf. Proc., \textbf{6}, Bar-Ilan Univ., Ramat Gan, 1993.%
\newline
\noindent \lbrack JudSh] H. Judah, S. Shelah, Baire property and axiom of
choice. \textsl{Israel J. Math.} \textbf{84} (1993), no. 3, 435--450.\newline
\noindent \lbrack JudSp] H. Judah, O. Spinas, Large cardinals and projective
sets. \textsl{Arch. Math. Logic} \textbf{36} (1997), no. 2, 137--155. 
\newline
\noindent \lbrack Kan] A. Kanamori, \textsl{The higher infinity. Large
cardinals in set theory from their beginnings}, Springer, 2$^{\text{nd}}$
ed. 2003 (1$^{\text{st}}$ ed. 1994).\newline
\noindent \lbrack KanM] A. Kanamori, M. Magidor, The evolution of large
cardinal axioms in set theory. \textsl{Higher set theory} (Proc. Conf.,
Math. Forschungsinst., Oberwolfach, 1977), pp. 99--275, Lecture Notes in
Math. \textbf{669}, Springer, 1978.\newline
\noindent \lbrack Kec1] A. S. Kechris, The axiom of determinacy implies
dependent choices in $L(\mathbb{R})$. \textsl{J. Symbolic Logic} \textbf{49.1%
} (1984), 161--173.\newline
\noindent \lbrack Kec2] A. S. Kechris: \textsl{Classical Descriptive Set
Theory.} Grad. Texts in Math. \textbf{156}, Springer, 1995.\newline
\noindent \lbrack KecS] A. S. Kechris, R. M. Solovay, On the relative
consistency strength of determinacy hypotheses. \textsl{Trans. Amer. Math.
Soc.} \textbf{290.1} (1985), 179--211.\newline
\noindent \lbrack Kei] H. J. Keisler, \textsl{Foundations of infinitesimal
calculus}, Prindle Weber and Schmidt, 1976.\newline
\noindent \lbrack Kel] J. L Kelley, \textsl{General topology.} Van Nostrand.
1955.\newline
\noindent \lbrack Kle] E. M. Kleinberg, \textsl{Infinitary combinatorics and
the axiom of determinateness.} Lecture Notes in Math. \textbf{612},
Springer, 1977.\newline
\noindent \lbrack KoeW] P. Koellner, W.H. Woodin, \textsl{Large cardinals
from determinacy}, Ch. 23 in [ForK].\newline
\noindent \lbrack Kun1] K. Kunen, Elementary embeddings and infinitary
combinatorics. \textsl{J. Symbolic Logic} \textbf{36} (1971), 407--413.%
\newline
\noindent \lbrack Kun2] K. Kunen, \textsl{Set theory. An introduction to
independence proofs. }Reprint of the 1980 original. Studies in Logic and the
Foundations of Mathematics, 102. North-Holland, 1983. \newline
\noindent \lbrack Kun3] K. Kunen, Random and Cohen reals. \textsl{Handbook
of set-theoretic topology}, (eds K. Kunen and J.E. Vaughan), 887--911,
North-Holland, 1984.\newline
\noindent \lbrack Kur] K. Kuratowski, \textsl{Topology}. Vol. I. (tr. J.
Jaworowski), Academic Press, PWN, 1966.\newline
\noindent \lbrack Lar] P. B. Larson, A brief history of determinacy, \textsl{%
The Cabal Seminar Vol. 4} (eds. A. S. Kechris, B. L\"{o}we, J. R. Steel),
Assoc. Symbolic Logic, 2010.\newline
\noindent \lbrack \L o\'{s}] J. \L o\'{s}, Quelques remarques, th\'{e}or\`{e}%
mes et probl\`{e}mes sur les classes d\'{e}finissables d'alg\`{e}bres, pp.
98-113 in \textsl{Mathematical Interpretations of Formal Systems},
North-Holland, 1955.\newline
\noindent \lbrack Lus] N. N. Lusin, \textsl{Le\c{c}ons sur les ensembles
analytiques}. Gauthier-Villars, 1930. \newline
\noindent \lbrack LusS] N. N. Lusin, W. Sierpi\'{n}ski, Sur quelques propri%
\'{e}t\'{e}s des ensembles (A), \textsl{Bull. Acad. Sci. Crac., Sc.Math. Nat.%
},\textsl{\ S\'{e}r A}, (1918) 35-48.\newline
\noindent \lbrack Mar] D. A. Martin, Measurable cardinals and analytic
games. \textsl{Fund. Math.} \textbf{66} (1969/1970) 287--291.\newline
\noindent \lbrack MarK] D. A. Martin, A. S. Kechris, Infinite games and
effective descriptive set theory, in: [Rog, part 4]\newline
\noindent \lbrack MarS] D. A. Martin, R. M. Solovay, Internal Cohen
extensions. \textsl{Ann. Math. Logic} \textbf{2.2} (1970), 143--178.\newline
\noindent \lbrack MarSt] D. A. Martin, J. R. Steel, A proof of projective
determinacy, \textsl{J. Amer. Math. Soc.} \textbf{2 }(1989) 71-125.\newline
\noindent \lbrack Mat] A. R. D. Mathias, Surrealist landscape with figures
(a survey of recent results in set theory), \textsl{Period. Math. Hung.} 
\textbf{10} (1979), 109-175.\newline
\noindent \lbrack Mau] R. D. Mauldin, ed., \textsl{The Scottish Book}, Birk%
\"{a}user, Boston, 1981.\newline
\noindent \lbrack vMil] J. van Mill, A note on the Effros Theorem, \textsl{%
Amer. Math. Monthly} \textbf{111.9} (2004), 801-806.\newline
\noindent \lbrack MilO] H. I. Miller and A. J. Ostaszewski, Group action and
shift-compactness, \textsl{J. Math. Anal. App.} \textbf{392} (2012), 23-39.%
\newline
\noindent \lbrack MonZ] D. Montgomery, L. Zippin, \textsl{Topological
transformation groups.} Krieger, 1974 (1$^{\text{st}}$ printing:
Interscience, 1955.)\newline
\noindent \lbrack Moo] J. Tatch Moore, The proper forcing axiom. \textsl{%
Proc. ICM Vol II}, 3--29, Hindustan Book Agency, 2010.\newline
\noindent \lbrack Mor] M. Morley, Homogeneous sets, 181-196, in: \textsl{%
Handbook of Mathematical Logic}, North-Holland, 1977.\newline
\noindent \lbrack Mos1] Y. N. Moschovakis, Uniformization in a playful
universe. \textsl{Bull. Amer. Math. Soc.} \textbf{77} (1971) 731--736.%
\newline
\noindent \lbrack Mos2] Y. N. Moschovakis, \textsl{Notes on set theory.} 2$^{%
\text{nd}}$ ed. Undergrad. Texts in Math., Springer, 2006.\newline
\noindent \lbrack Most] A. Mostowski, On the principle of dependent choices. 
\textsl{Fund. Math.} \textbf{35} (1948). 127--130.\newline
\noindent \lbrack Myc] J. Mycielski, On the axiom of determinateness. 
\textsl{Fund. Math.} \textbf{53} (1963/19640, 205--224.\newline
\noindent \lbrack MycS] J. Mycielski, H. Steinhaus, A mathematical axiom
contradicting the axiom of choice. \textsl{Bull. Acad. Polon. Sci. S\'{e}r.
Sci. Math. Astronom. Phys.} \textbf{10} (1962), 1--3.\newline
\noindent \lbrack MycSw] J. Mycielski, S. \'{S}wierczkowski, On the Lebesgue
measurability and the axiom of determinateness. \textsl{Fund. Math.} \textbf{%
54} (1964), 67--71.\newline
\noindent \lbrack MycT] J. Mycielski and G. Tomkowicz, Shadows of the axiom
of choice in the universe $L(\mathbb{R}).$ \textsl{Arch. Math. Logic} (doi
10.1007/s00153-017-0596-x).\newline
\noindent \lbrack MyhS] J. Myhill, D. Scott, Ordinal definability, in: 
\textsl{Axiomatic Set Theory}, Proc Symp. Pure Math Vol XIII-I, Amer. Math.
Soc. (1971), 271-278.\newline
\noindent \lbrack Nee] I. Neeman, \textsl{Determinacy in }$L(\mathbb{R})$,
Ch. 21 in [ForK].\newline
\noindent \lbrack Neu1] J. von Neumann, \"{U}ber die Definition durch
transfinite Induktion und verwandte Fragen der allgemeinen Mengenlehre. 
\textsl{Math. Annalen} \textbf{99} (1928), 373-391 (\textsl{Collected Works,
Vol. I} (ed. A. H. Taub), Pergamon Press, 1961, 320-338).\newline
\noindent \lbrack Neu2] J. von Neumann, Der Axiomatisierung der Mengenlehre. 
\textsl{Math. Z.} \textbf{27}, (1928), 669-752 (\textsl{Works I}, 339-423).%
\newline
\noindent \lbrack Nik] O. Nikodym, Sur une propri\'{e}t\'{e} de l'op\'{e}%
ration A. \textsl{Fund. Math.} \textbf{7} (1925), 149-154.\newline
\noindent \lbrack Ost1] A. J. Ostaszewski, On countably compact, perfectly
normal spaces. \textsl{J. London Math. Soc.} (2) \textbf{14} (1976),
505--516.\newline
\noindent \lbrack Ost2] A. J. Ostaszewski, On how to trap a gap: "An
introduction to independence for Analysts by H.G. Dales and W.H. Woodin", 
\textsl{Bull. London Math. Soc.} \textbf{21 }(1989), 197-208 (Extended
Review Article of [DalW1]).\newline
\noindent \lbrack Ost3] A. J. Ostaszewski, Analytic Baire spaces. \textsl{%
Fund. Math.} \textbf{217} (2012), no. 3, 189--210.\newline
\noindent \lbrack Ost4] A. J. Ostaszewski, Almost completeness and the
Effros Theorem in normed groups, \textsl{Topology Proceedings} \textbf{41}
(2013), 99-110 (fuller version: arXiv.1606.04496).\newline
\noindent \lbrack Ost5] A. J. Ostaszewski, Shift-compactness in almost
analytic submetrizable Baire groups and spaces, survey article, \textsl{%
Topology Proceedings} \textbf{41} (2013), 123-151.\newline
\noindent \lbrack Ost6] A. J. Ostaszewski, Effros, Baire, Steinhaus and
non-separability, \textsl{Topology and its Applications}, Mary Ellen Rudin
Special Issue, \textbf{195 }(2015), 265-274.\newline
\noindent \lbrack Ost7] A. J. Ostaszewski, Topological descriptive set
theory, in: [FalGOS], 426-429.\newline
\noindent \lbrack Oxt] J. C. Oxtoby: \textsl{Measure and category}, 2nd ed.
Graduate Texts in Math. \textbf{2}, Springer, 1980 (1$^{\text{st}}$ ed.
1972).\newline
\noindent \lbrack Paw] J. Pawlikowski, Lebesgue measurability implies Baire
property. \textsl{Bull. Sci. Math.} (2) \textbf{109} (1985), 321--324\newline
\noindent \lbrack Pin1] D. Pincus, Independence of the prime ideal theorem
from the Hahn Banach theorem. \textsl{Bull. Amer. Math. Soc.} \textbf{78}
(1972), 766--770.\newline
\noindent \lbrack Pin2] D. Pincus, The strength of the Hahn-Banach theorem.%
\textsl{\ Victoria Symposium on Nonstandard Analysis }(Univ. Victoria,
Victoria, B.C., 1972), pp. 203--248, Lecture Notes in Math. \textbf{369},
Springer 1974.\newline
\noindent \lbrack Pin3] D. Pincus, Adding dependent choice to the prime
ideal theorem. Logic Colloquium 76, \textsl{Studies in Logic and Found. Math.%
} Vol. \textbf{87}, 547--565, North-Holland, 1977.\newline
\noindent \lbrack Pin4] D. Pincus, Adding dependent choice. \textsl{Ann.
Math. Logic} \textbf{11} (1977), no. 1, 105--145.\newline
\noindent \lbrack PinS] D. Pincus, R. Solovay, Definability of measures and
ultrafilters.\textsl{\ J. Symbolic Logic} \textbf{42.2} (1977), 179--190.%
\newline
\noindent \lbrack Rai] J. Raisonnier, A mathematical proof of S. Shelah's
theorem on the measure problem and related results. \textsl{Israel J. Math.} 
\textbf{48} (1984), no. 1, 48--56.\newline
\noindent \lbrack RaiS] J. Raisonnier, J. Stern, Mesurabilit\'{e} and propri%
\'{e}t\'{e} de Baire, \textsl{Comptes Rendus Acad. Sci. I. (Math.)} \textbf{%
296} (1983), 323-326.\newline
\noindent \lbrack Ram] F. P. Ramsey, On a problem of formal logic, \textsl{%
Proc. London Math. Soc.}, \textbf{30} (1929), 338-384.\newline
\noindent \lbrack Rob1], A. Robinson, \textsl{Introduction to model theory
and to the metamathematics of algebra}, North-Holland, 1965 (1$^{\text{st}}$
ed 1963).\newline
\noindent \lbrack Rob2], A. Robinson, \textsl{Non-standard analysis},
North-Holland, 1970.\newline
\noindent \lbrack Rog] C. A. Rogers et al., \textsl{Analytic sets}, Academic
Press, 1980.\newline
\noindent \lbrack Rud] W. Rudin, \textsl{Functional analysis.} 2$^{\text{nd}%
} $ ed. International Series in Pure and Applied Mathematics. McGraw-Hill,
1991 (1$^{\text{st}}$ ed. 1973).\newline
\noindent \lbrack Sco1] D. Scott, measurable cardinals and constructible
sets, Bull. Acad. Polon. Sci., \textbf{7} (1961), 145-149.\newline
\noindent \lbrack Sco2] D. Scott, Axiomatizing set theory, in: \textsl{%
Axiomatic Set Theory}, Proc Symp. Pure Math Vol XIII-II, Amer. Math. Soc.
(1974), 207-214.\newline
\noindent \lbrack She1] S. Shelah, Can you take Solovay's inaccessible away? 
\textsl{Israel J. Math.} \textbf{48} (1984), 1-47.\newline
\noindent \lbrack She2] S. Shelah, On measure and category. \textsl{Israel
J. Math.} \textbf{52} (1985), no. 1-2, 110--114.\newline
\noindent \lbrack SheW] S. Shelah, H. Woodin, Large cardinals imply that
every reasonably definable set of reals is Lebesgue measurable. \textsl{%
Israel J. Math.} \textbf{70.3} (1990), 381--394.\newline
\noindent \lbrack Sie] W. Sierpi\'{n}ski, \textsl{Hypoth\`{e}se du continu},
2$^{\text{nd}}$ ed., Chelsea, 1956\newline
\noindent \lbrack Sil] J. H. Silver, Some applications of model theory in
set theory. \textsl{Ann. Math. Logic} \textbf{3} (1971), 45--110.\newline
\noindent \lbrack Sko] T. Skolem, \"{U}ber die Nicht-Charakterisierbarkeit
der Zahlenreihe mittels endlich oder abz\"{a}hlbar unendlich vieler Aussagen
mit ausschliesslich Zahlenvariablen, \textsl{Fund. Math.} \textbf{23}
(1934), 150-161.\newline
\noindent \lbrack Sol1] R. M. Solovay, $2^{\aleph _{0}}$ can be anything it
ought to be, \textsl{The theory of models}, 435 (Proc. 1963 Int. Symp.
Berkeley) Studies in Logic and the Foundations of Mathematics, eds. J. W.
Addison, Leon Henkin, Alfred Tarski, North-Holland, 1965.\newline
\noindent \lbrack Sol2] R. M. Solovay, On the cardinality of $\Sigma
_{2}^{1} $ sets of reals. \textsl{Foundations of Mathematics} (Symposium
Commemorating Kurt G\"{o}del, Columbus, Ohio, 1966) pp. 58--73, Springer,
1969.\newline
\noindent \lbrack Sol3] R. M. Solovay, A model of set-theory in which every
set of reals is Lebesgue measurable. \textsl{Ann. of Math.} (2) \textbf{92}
(1970), 1--56.\newline
\noindent \lbrack Sol4] R. M. Solovay, The independence of DC from AD. 
\textsl{The Cabal Seminar 76--77} (Proc. Caltech-UCLA Logic Sem., 1976--77),
pp. 171--183, Lecture Notes in Math. \textbf{689}, Springer, 1978.\newline
\noindent \lbrack SolAH] Solovay, Robert M.; Arthan, R. D.; Harrison, John
Some new results on decidability for elementary algebra and geometry. Ann.
Pure Appl. Logic 163 (2012), no. 12, 1765--1802.\newline
\noindent \lbrack SolRK] R. M. Solovay, W. N. Reinhardt, A. Kanamori, Strong
axioms of infinity and elementary embeddings, \textsl{Ann. Math. Logic} 
\textbf{13} (1978), 73-116.\newline
\noindent \lbrack SolT] R. M. Solovay, S. Tennenbaum, Iterated Cohen
extensions and Souslin's problem. \textsl{Ann. of Math.} (2) \textbf{94}
(1971), 201--245.\newline
\noindent \lbrack Sou] M. Souslin, Sur une d\'{e}finition des ensembles
mesurables B sans nombres transfinis, \textsl{C.R Acad. Sci.} \textbf{164}
(1917) 88-91. \newline
\noindent \lbrack Ste] J. R. Steel, G\"{o}del's program. \textsl{%
Interpreting G\"{o}del}, 153--179, Cambridge Univ. Press, 2014.\newline
\noindent \lbrack Ster] J. Stern, Regularity properties of definable sets of
reals, \textsl{Annals of Pure and Appl. Logic}, \textbf{29} (1985), 289-324.%
\newline
\noindent \lbrack Sto] A. H. Stone, Analytic sets in non-separable metric
spaces, [Rog, Part 5], 471-480.\newline
\noindent \lbrack Tar1] A. Tarski, Une contribution \`{a} la th\'{e}orie de
la mesure, \textsl{Fund. Math.} \textbf{15} (1930), 42-50.\newline
\noindent \lbrack Tar2] A. Tarski, Axiomatic and algebraic aspects of two
theorems on sums of cardinals. \textsl{Fund. Math.} \textbf{35} (1948),
79--104.\newline
\noindent \lbrack Tod] S. Todor\v{c}evi\'{c}. Generic absoluteness and the
continuum. \textsl{Mathematical Research Letters}, 9:465--472, 2002.\newline
\noindent \lbrack TomW] G. Tomkowicz and S. Wagon, \textsl{The Banach-Tarski
paradox}, 2$^{\text{nd}}$ ed., Cambridge University Press, 2016 (1st ed., S.
Wagon, CUP, 1985).\newline
\noindent \lbrack Ula] S. Ulam, Zur Masstheorie in der allgemeinen
Mengenlehre. \textsl{Fund. Math.} \textbf{16} (1930), 140-150.\newline
\noindent \lbrack Vau] R. L. Vaught, Alfred Tarski's work in model theory, 
\textsl{J. Symbolic Logic}, \textbf{51} (1986), 869-882.\newline
\noindent \lbrack Vel] B. Veli\v{c}kovi\'{c}, Forcing axioms and stationary
sets. \textsl{Adv. Math.} \textbf{94} (2) (1992), 256-- 284\newline
\noindent \lbrack Wei] W. Weiss, Versions of Martins' Axiom. \textsl{%
Handbook of set-theoretic topology} (eds K. Kunen and J.E. Vaughan),
827--886, North-Holland, 1984.\newline
\noindent \lbrack Wol] E. Wolk, On the principle of dependent choices and
some forms of Zorn's lemma, \textsl{Canad. Bull. Math.}, \textbf{26} (1983),
365-367.\newline
\noindent \lbrack Woo1] W. H. Woodin, \textsl{The Axiom of Determinacy,
Forcing Axioms, and the Nonstationary Ideal}, de Gruyter Series in Logic and
its Applications 1, de Gruyter, 1999.\newline
\noindent \lbrack Woo2] W. H. Woodin, The continuum hypothesis. I; II, 
\textsl{Notices Amer. Math. Soc. }\textbf{48} (2001), 567--576; 681--690.%
\newline
\noindent \lbrack Woo3] W. H. Woodin, The weak ultimate L conjecture. 
\textsl{Infinity, computability, and metamathematics}, 309--329, Tributes,
23, Coll. Publ., London, 2014. (MR3307892)\newline
\noindent \lbrack Wri] J. D. Maitland Wright, Functional Analysis for the
practical man, 283--290 in \textsl{Functional Analysis: Surveys and Recent
Results} \textbf{27}, North-Holland Math. Studies, 1977.\newline
\bigskip

\bigskip

\noindent Mathematics Department, Imperial College, London SW7 2AZ;
n.bingham@ic.ac.uk \newline
Mathematics Department, London School of Economics, Houghton Street, London
WC2A 2AE; A.J.Ostaszewski@lse.ac.uk

\end{document}